\documentclass[12pt]{article}

\usepackage{amsfonts,amssymb,amsthm,graphicx,rotating,appendix,xcolor,helvet,amsmath}
\usepackage{bm}
\usepackage[utf8]{inputenc}
\usepackage[english]{babel}
\usepackage{algorithm}
\floatname{algorithm}{Algorithm}
\usepackage{booktabs}
\usepackage{rotating} 
\usepackage{algpseudocode}
\usepackage{caption}

\usepackage[all,cmtip]{xy}

\usepackage[round,authoryear]{natbib}

\usepackage{color}
\usepackage{xcolor}
\usepackage{multirow}
\usepackage{changepage}

\usepackage{todonotes}

\usepackage{lineno}
\usepackage{booktabs}

\DeclareMathOperator{\mcF}{\mathcal{F}}

\DeclareMathOperator{\mbR}{\mathbb{R}}
\def\D{\mathrm{d}}

\newtheorem{prop}{Proposition}

\newtheorem{lemma}{Lemma}

\usepackage{enumitem}
\usepackage[backref=page,
            colorlinks,
            citecolor=blue,
            linkcolor=blue,
            pdfborder={0 0 0}]{hyperref} 
\usepackage{orcidlink}
\numberwithin{equation}{subsection}

\begin{document}

\title{Inference for two-stage sampling in spatial surveys}

\author{
Guillaume Chauvet \and Olivier Bouriaud \and Trinh H.K. Duong
}


\maketitle

\abstract{This paper develops a design-based asymptotic theory for two-stage sampling over a continuous spatial domain. The target parameters are integral totals, or smooth functions of such totals, defined over a fixed bounded territory partitioned into an increasingly fine collection of primary sampling units. Within this fixed-area framework, we derive the order of the variance components of the Horvitz–Thompson estimator when secondary sampling units are points selected from continuous sub-regions. We establish design consistency of the estimator, as well as consistency of variance estimators under explicit regularity conditions on inclusion probabilities, sampling densities, and the study variable. The results are extended to plug-in estimators of smooth scale-invariant functions of totals. For high-entropy first-stage designs, we further show that a H\'{a}jek-type variance estimator based only on first-order inclusion probabilities is consistent, and we prove asymptotic normality of both total and plug-in estimators. A simulation study inspired by forest inventory applications illustrates the finite-sample performance of the proposed estimators and variance estimators.}


\section{Introduction} \label{sec:intro}

\noindent Multistage sampling is a central component of design-based inference for complex surveys. In the classical finite-population framework, primary sampling units (PSUs) are selected at a first stage and secondary sampling units (SSUs) are drawn conditionally at subsequent stages \citep{sarswewre92,loh:10}. For such designs, the asymptotic theory is well established \citep{kre:rao:81,ohl:89}. More recently, \citet{chauvet2020inference} established general high-level conditions ensuring consistency and asymptotic normality of Horvitz–Thompson (HT) estimators under multistage sampling. \\
\noindent These results are formulated for discrete finite populations at all stages. In many applications, however, the target of inference is defined over a continuous spatial domain. Examples include forest inventories, soil surveys, and environmental monitoring, where a territory is partitioned into spatial sub-regions at the first stage (for example, grid cells), and sampling locations are selected within sampled regions at the second stage. Field measurements performed at sampled locations are used to estimate integral parameters of interest over the territory, such as total biomass, total volume, or proportions of area. A design-based approach for such spatial settings can be traced back to extensions of the HT principle to continuous universes \citep{cor:93}, and to subsequent developments in spatial and environmental statistics \citep{fattorini2004variance,fattorini2011two,fattorini2017design}. \\
\noindent Two-stage sampling when the PSUs are continuous sub-regions has been considered in \cite{fattorini2004variance}. These results are extended in \cite{duong24}, who developed a framework for two-stage sampling schemes, where two-phase sampling can be applied inside the sampled sub-regions. Despite this literature, a general asymptotic theory for two-stage sampling over continuous territories, parallel to that available for discrete multistage designs, remains incomplete. In the continous setting, second-stage inclusion probabilities are replaced by inclusion density functions and totals are expressed as integrals. Establishing consistency and asymptotic normality of estimators therefore requires explicit control of variability at both stages under assumptions adapted to continous domains. The present paper complements and extends \cite{duong24} by providing a formal design-based asymptotic theory for two-stage spatial sampling designs, including consistency, variance estimation and asymptotic normality. Compared with \cite{chauvet2020inference}, the present work does not merely transfer two-stage finite-population asymptotics to a continuous spatial setting. It also develops a design-based theory for plug-in estimators of smooth scale-invariant functions of spatial totals, including linearization and variance estimation, which was not addressed in the earlier finite-population framework.\\
\noindent The objective of this paper is to develop a design-based asymptotic theory for two-stage sampling over a continuous territory. We adopt the fixed-area asymptotic framework of \citet{fattorini2017design}, in which the domain is partitioned into an increasingly fine collection of PSUs while its total measure remains fixed. Within this setting, we establish the following results. First, we derive a decomposition of the variance of the two-stage HT estimator into first-stage and second-stage components adapted to continuous secondary units, and determine their respective orders of magnitude under mild regularity conditions. This clarifies the relative contribution of the two stages, in particular when the number of within-PSU sampling locations is bounded. Secondly, we prove $\sqrt{n_I}$-consistency of the HT estimator for integral totals, together with consistency of the associated HT and Yates–Grundy (YG) variance estimators under explicit assumptions on inclusion probabilities and sampling densities. We focus on two-stage designs in which the within-PSU sampling design allows the construction of second-order inclusion densities. The practically important one-point-per-cell case, as well as post-stratified and more general two-phase extensions, are not covered by the present asymptotic theory and are left for future work. Thirdly, we extend these results to smooth scale-invariant functions of totals. Under suitable regularity conditions on the functional of interest, we establish linearization-based approximations and consistency of the corresponding variance estimators. Finally, for high-entropy first-stage designs, we show that the first-stage variance component can be consistently estimated using a Hájek-type estimator that depends only on first-order inclusion probabilities. In this setting, we establish asymptotic normality of HT and plug-in estimators, thereby providing a rigorous justification for normal-approximation confidence intervals.

\section{Notation} \label{sec:notation}

\noindent We consider a bounded measurable set $\mcF \subset \mbR^d$, endowed with the Lebesgue measure, whose total measure is denoted by $A$. The territory is partitioned into $N_I$ disjoint measurable sub-regions, forming a finite population $U_I=\{u_1,\ldots,u_{N_I}\}$ of PSUs. They need not have identical shapes or areas. We denote by $A_i$ the measure of $u_i$. \\
\noindent To study asymptotic properties, we adopt the fixed-area consistency framework of \citet{fattorini2017design}, under which the total measure $A$ is fixed while the number of PSUs satisfies $N_I \to \infty$. Consequently, the partition becomes increasingly fine and PSU diameters shrink to zero. This regime is natural in spatial sampling, and all asymptotic results are derived under this framework.

\subsection{Estimation of a total} \label{sec:notation:ssec1}

\noindent Let $\rho: \mcF \to \mbR$ be a Lebesgue integrable function. The target parameter is the integral
    \begin{equation} \label{sec:notation:ssec1:eq1}
    \tau_{\rho} = \int_{\mcF} \rho(x) \, \D{x}.
    \end{equation}
This formulation encompasses both intrinsic spatial variables and totals over finite populations $U$ embedded in $\mcF$. In the latter case, following \cite{chaboubri23} and \cite{bouriaud2024wsm}, we represent unit-level attributes $y_k,~k \in U$, through a non-negative link function $L(\cdot,\cdot):\mathcal{F} \times U \rightarrow \mathbb{R}^+$ and define
    \begin{align} \label{weight:share}
    \rho(x) = \sum_{k \in U} \dfrac{L(x,k) y_k}{M_{+k}}, & &
    M_{+k} = \int_{\mcF} L(x,k) \, \D{x},
    \end{align}
so that $\tau_{\rho} = \sum_{k \in U} y_k$ provided that $M_{+k}>0$ for all $k \in U$. \\
\noindent At the first-stage, a sample $S_I \subset U_I$ of size $n_I$ is selected according to a sampling design $p_I(\cdot)$, with inclusion probabilities $\pi_{Ii}$ and joint inclusion probabilities $\pi_{Iij}$.
Conditionally on $S_I$, for each selected PSU $u_i$, a second-stage sample $S_{IIi}$ of $n_{IIi}$ locations is drawn independently from a density $f_{IIi}$ supported on $u_i$. The associated conditional inclusion density function is $\pi_{IIi}(x) = n_{IIi} f_{IIi}(x)$. The overall sample is $S = \bigcup_{u_i \in S_I} S_{IIi}$. We assume invariance of the second-stage design, meaning that the second-stage sampling mechanism does not depend on the realization of $S_I$ \citep[][chapter 4]{sarswewre92}. In addition, the second-stage samples $S_{IIi}$ are assumed to be conditionally independent across PSUs, given $S_I$.\\
\noindent The associated HT estimator of $\tau_{\rho}$ \citep{cor:93,fattorini2004variance,duong24} is
    \begin{eqnarray} \label{sec:notation:ssec1:eq8}
        \Hat{\tau}_{\rho} & = & \sum_{u_i \in S_I} \frac{\Hat{\tau}_{\rho i}}{\pi_{Ii}}
        ~\textrm{  with  }~ \Hat{\tau}_{\rho i} = \sum_{x \in S_{IIi}} \frac{\rho(x)}{\pi_{IIi}(x)},
    \end{eqnarray}
where $\Hat{\tau}_{\rho i}$ is an estimator of the sub-total $\tau_{\rho i} = \int_{u_i} \rho(x) \, \D{x}$. The variance of $\Hat{\tau}_{\rho}$ admits the decomposition
    \begin{align} \label{sec:notation:ssec1:eq9}
    V(\Hat{\tau}_{\rho}) & = V_1(\Hat{\tau}_{\rho})+V_2(\Hat{\tau}_{\rho})+V_3(\Hat{\tau}_{\rho}),
    \end{align}
where $V_1(\Hat{\tau}_{\rho})$ corresponds the first-stage variability and $V_2(\Hat{\tau}_{\rho})+V_3(\Hat{\tau}_{\rho})$ to second-stage contributions; explicit expressions for these quantities are recalled in Appendix \ref{app:add:not:1}. Controlling the relative magnitude of these components is important for the asymptotic analysis developed below. \\
\noindent We consider two variance estimators corresponding to the decomposition of $V(\Hat{\tau}_{\rho})$. The HT variance estimator is defined as
    \begin{align} \label{sec:notation:ssec1:eq12}
    \Hat{V}_{HT}\left(\Hat{\tau}_{\rho}\right) & = \Hat{V}_{HT,A}\left(\Hat{\tau}_{\rho}\right) + \Hat{V}_{B}\left(\Hat{\tau}_{\rho}\right),
    \end{align}
where $\Hat{V}_{HT,A}\left(\Hat{\tau}_{\rho}\right)$ estimates the first-stage contribution and $\Hat{V}_{B}\left(\Hat{\tau}_{\rho}\right)$ estimates the second-stage component. For fixed-size sampling designs at the first-stage, The YG variance estimator may alternatively be used. It is defined by
    \begin{align} \label{sec:notation:ssec1:eq14}
    \Hat{V}_{YG}\left(\Hat{\tau}_{\rho}\right)
    & = \Hat{V}_{YG,A}\left(\Hat{\tau}_{\rho}\right) + \Hat{V}_{B}\left(\Hat{\tau}_{\rho}\right),
    \end{align}
where $\Hat{V}_{YG,A}\left(\Hat{\tau}_{\rho}\right)$ replaces the first-stage HT component in (\ref{sec:notation:ssec1:eq12}) by the corresponding YG form. Explicit expressions of these estimators are given in Appendix~\ref{app:add:not:2}.

\subsection{Plug-in estimation} \label{sec:notation:ssec2}

\noindent We consider parameters of the form
    \begin{equation*}
    \theta=g(\bm{\tau}_{\bm{\rho}}),
    \end{equation*}
where $\bm{\rho}(x) \in \mbR^Q$ is a vector of measurable functions and $g:\mathbb{R}^Q \to \mathbb{R}$ is differentiable. The natural estimator is the plug-in estimator
$\hat{\theta} = g(\hat{\bm{\tau}}_{\bm{\rho}})$, where $\hat{\bm{\tau}}_{\bm{\rho}}$ denotes the vector of HT estimators. The linearized variable associated with $\theta$ is
    \begin{eqnarray} \label{sec:notation:ssec2:eq1}
    l(x) & = & \{\bm{g}'(\bm{\tau}_{\bm{\rho}})\}^{\top} \bm{\rho}(x),
    \end{eqnarray}
and the corresponding linearization-based variance approximation is obtained by replacing $\rho(x)$ with $l(x)$ in the variance decomposition (\ref{sec:notation:ssec1:eq9}). The associated variance estimators are defined by replacing $\rho(x)$ with
    \begin{eqnarray} \label{sec:notation:ssec2:eq2}
    \hat{l}(x) & = & \{\bm{g}'(\hat{\bm{\tau}}_{\bm{\rho}})\}^{\top} \bm{\rho}(x).
    \end{eqnarray}
in the HT or in the YG variance estimator.

\section{Assumptions} \label{sec:assump}

\noindent To derive consistency and asymptotic normality results, we impose regularity conditions on the first-stage design, the second-stage design, and the variable of interest. These conditions ensure control of higher-order inclusion moments, regularity of spatial sampling densities, and stability of the variance component. We denote by
    \begin{align} \label{sec:assump:eq1}
        A_0 = \frac{1}{N_I} \sum_{i=1}^{N_I} A_i = \frac{A}{N_I}, & \quad
        n_{II0} = \frac{1}{N_I} \sum_{i=1}^{N_I} n_{IIi},
    \end{align}
the average share of the territory measure per PSU and the average number of second-stage points selected per PSU. Under the fixed-area framework \citep{fattorini2017design}, we assume that $N_I \to \infty$ and $A_0 \to 0$ as the partition of the territory becomes increasingly fine, and that the number of sampled PSUs satisfied $n_I \to \infty$.

\subsection{First-stage sampling design}

\noindent The following conditions extend those of \citet{chauvet2020inference} to accommodate the continuous-domain structure and the control of higher-order moments required for plug-in estimation.
    \begin{itemize}
      \item[FS1:] The first-stage sampling fraction is bounded away from one, and inclusion probabilities are of order $n_I/N_I$; that is, there exists constants $f_{I0}<1$ and $0<c_{I1}\leq C_{I1}< \infty$ such that
        \begin{equation} \label{sec:assump:eq2}
          \frac{n_I}{N_I} \leq f_{I0}, \qquad  c_{I1} \leq \frac{N_I}{n_I} \pi_{Ii} \leq C_{I1}.
        \end{equation}
      \item[FS2:] Higher-order centered inclusion moments up to order four are uniformly bounded at rates consistent with simple random or high-entropy sampling. More precisely, there exist constants $C_{I2},C_{I3},C_{I4}>0$ such that the following bounds hold:
        \begin{align}
          \Delta_{I2} & \equiv \max_{i \neq {j}=1,\ldots,N_I} \left| E\left\{(I_{Ii}-\pi_{Ii})(I_{Ij}-\pi_{Ij})\right\} \right| \leq C_{I2} N_I^{-2} n_I, \nonumber \\
          \Delta_{I3} & \equiv \max_{i \neq {j} \neq {i'} =1,\ldots,N_I} \left| E\left\{(I_{Ii}-\pi_{Ii}) (I_{Ij}-\pi_{Ij})(I_{Ii'}-\pi_{Ii'})\right\} \right| \leq C_{I3} N_I^{-3} n_I^2, \label{sec:assump:eq4} \\
          \Delta_{I4} & \equiv \max_{i \neq {j} \neq {i'} \neq {j'}=1,\ldots,N_I} \left| E\left\{(I_{Ii}-\pi_{Ii})(I_{Ij}-\pi_{Ij})(I_{Ii'}-\pi_{Ii'})(I_{Ij'}-\pi_{Ij'})\right\} \right| \leq C_{I4} N_I^{-4} n_I^2, \nonumber
        \end{align}
      where $I_{Ii}$ denote the sample membership indicator of PSU $u_i$ in $S_I$.
      \item[FS3:] There exists $c_{I2}>0$ such that for any $i \ne j =1,\ldots,N_I$:
        \begin{eqnarray} \label{sec:assump:eq5}
           c_{I2} N_I^{-2} n_I^2 & \leq \pi_{Iij}.
        \end{eqnarray}
    \end{itemize}

\noindent Assumption (FS1) was done in \citet{chauvet2020inference}. It concerns first-order inclusion probabilities and is controlled by the survey designer. Assumption (FS2) involves higher-order inclusion probabilities (orders 2-4) and holds for high-entropy designs such as simple random sampling or rejective sampling \citep{Haj64,BoiLopRui12}. Conditions (FS1)–(FS2) imply Assumption 2 of \citet{chauvet2020inference}, but additionally ensure control of fourth-order moments of the HT estimator. This strengthening is required to establish mean-square consistency of variance estimators and of plug-in functionals. Assumption (FS3) prevents instability of variance estimators; see \citet[][Section 3.1]{chauvet2020inference} for a more detailed discussion.

\subsection{Second-stage sampling design}

\noindent To control within-PSU variability under the fixed-area framework, we impose uniform regularity conditions on second-stage sample sizes, PSU areas, and sampling densities.
    \begin{itemize}
      \item[SS1:] There exist some constants $0<\lambda_1 \leq \Lambda_1 < \infty$, $0<\phi_1 \leq \Phi_1 \leq \infty$, and $0<c_1 \leq C_1 \leq \infty$, such that for any PSU $u_i$:
        \begin{eqnarray}
        \lambda_1 n_{II0} & \leq n_{IIi} & \leq \Lambda_1 n_{II0},
        \label{sec:assump:eq6} \\
        \phi_1 \frac{A}{N_I} & \leq A_i & \leq \Phi_1 \frac{A}{N_I}, \label{sec:assump:eq7} \\
        \frac{c_1}{A_i} & \leq f_{IIi}(x) & \leq \frac{C_1}{A_i} \quad \text{for any } x \in u_i. \label{sec:assump:eq8}
        \end{eqnarray}
      Also, $n_{IIi} \geq 2$ for any PSU $u_i \in U_I$.
    \end{itemize}
These conditions ensure that second-stage sample sizes $n_{IIi}$ and PSU measures $A_i$ remain comparable across the partition. This holds automatically under regular grid partitions and when selecting the same number of points per PSU. The marginal density $f_{IIi}(\cdot)$ is assumed to be uniformly bounded from zero and infinity; in particular, the condition is satisfied when $f_{IIi}(x)=A_i^{-1}$, i.e., under independent uniform sampling.

\subsection{Assumptions on the variable and on the functional}

\noindent For total estimation, we require moment and non-degeneracy conditions.
    \begin{itemize}
      \item[VA1:] The function $\rho$ has a bounded fourth moment and a mean bounded away from zero; that is, there exist constants $m_1,M_1$ such that
         \begin{equation}
         \frac{1}{A} \int_{\mcF} \{\rho(x)\}^4 \D{x} \leq M_1, \quad
         \left|\frac{1}{A} \int_{\mcF} \rho(x) \D{x} \right| \geq m_1 >0.  \label{sec:assump:eq10}
         \end{equation}
      \item[VA2:] The leading first-stage variance component remains asymptotically non-negligible; that is, there exists some constant $m_2$ such that
         \begin{align} \label{sec:assump:eq11}
         \frac{n_I}{A^{2}} V_1(\Hat{\tau}_{\rho}) \geq m_2 >0.
         \end{align}
    \end{itemize}

\noindent For parameters of the form $\theta=g(\bm{\tau}_{\bm{\rho}})$, we impose smoothness conditions on $g$.
    \begin{itemize}
      \item[PI1:] The function $g(\cdot)$ is homogeneous of degree $\beta=0$; that is, $g(r\mathbf{a})=r^{\beta} g(\mathbf{a})= g(\mathbf{a})$ for any real $r>0$ and any $\mathbf{a} \in \mathbb{R}^Q$.
      \item[PI2:] The attribute density function $\bm{\rho}(\cdot)$ takes values in a convex compact set $\mathcal{C} \subset \mathbb{R}^Q$, and $\bm{\mu}_{\bm{\rho}} = A^{-1} \bm{\tau}_{\bm{\rho}}$ is an interior point of $\mathcal{C}$. The function $g(\cdot)$ is differentiable on $\mathcal{C}$, and the differential $\mathbf{g}'$ is locally Lipschitz on $\mathcal{C}$; that is, there exists $K_g>0$ such that
        \begin{equation*}
        \|\mathbf{g}'(\mathbf{b})-\mathbf{g}'(\mathbf{a})\| \leq K_g \|\mathbf{b}-\mathbf{a}\| \quad \text{for any } \mathbf{a},\mathbf{b} \in \mathcal{C}^2,    \end{equation*}
      where $\|\cdot\|$ stands for the Euclidean norm.
    \end{itemize}
These conditions ensure validity of linearization and control of the remainder term in the delta expansion.
The plug-in theory is aimed at the main nonlinear parameters used in spatial survey applications. Together with totals and linear combinations of totals, which are covered by the results of Section~\ref{sec:notation:ssec1}, Assumptions (PI1)-(PI2) cover scale-invariant smooth functionals such as ratios, proportions, domain means, relative changes and related composition indicators. These include the principal estimands considered in forest and environmental inventories.

\section{Consistency of estimators} \label{sec:consist}

\noindent We first control the order of magnitude of each variance component. This allows us to establish $\sqrt{n_I}$-consistency of the estimator. We then study consistency of variance estimators. Finally, we extend the results to plug-in functionals.

\subsection{Estimators of totals}


\begin{prop} \label{prop1}
Suppose that the first-stage assumptions (FS1)-(FS2), the second-stage assumption (SS1) and the assumption (VA1) hold. Then:
    \begin{align}
    V_1(\Hat{\tau}_{y}) & = O\left(\frac{A^2}{n_I}\right) \label{sec:consist:eq1} \\
    V_2(\Hat{\tau}_{y}) & = O\left(\frac{A^2}{n_I~n_{II0}}\right) \label{sec:consist:eq2} \\
    V_3(\Hat{\tau}_{y}) & = O\left(\frac{A^2}{N_I~n_{II0}}\right). \label{sec:consist:eq3}
    \end{align}
\end{prop}

\noindent The proof is provided in Appendix \ref{app:proof:prop1}. Proposition \ref{prop1} shows that the leading variance component is of order $A^2 n_{I}^{-1}$.
The second-stage contribution $V_2$ is reduced by a factor $n_{II0}$, while the term $V_3$ is further reduced by the factor $N_I^{-1}$. In particular; if $n_{II0} \to \infty$, the second-stage contribution becomes negligible; if $n_{II0}$ remains bounded, $V_1$ and $V_2$ are of the same order; if the first-stage sampling fraction satisfies $n_I/N_I \to 0$, the term $V_3$ is asymptotically negligible. \\
\noindent Proposition \ref{prop1} immediately implies the mean-square convergence of the HT estimator under the fixed-area asymptotic regime, as formalized in Proposition \ref{prop2}.

\begin{prop} \label{prop2}
  Suppose that assumptions (FS1)-(FS2), (SS1) and (VA1) hold. Then the HT estimator is design-unbiased. Also, we have
    \begin{eqnarray}
      E\left[ A^{-1} \left\{\Hat{\tau}_{y}-\tau_y\right\}\right]^2 = O(n_I^{-1})
      & \textrm{ and } & \frac{\Hat{\tau}_{y}}{\tau_y} \longrightarrow_{Pr} 1, \label{sec:consist:eq4}
    \end{eqnarray}
  where $\rightarrow_{Pr}$ stands for the convergence in probability.
\end{prop}

\noindent The $\sqrt{n_I}$-consistency of the HT variance estimator and of the YG variance estimator is established in Propositions \ref{prop3} and \ref{prop4}, respectively. The proof of Proposition \ref{prop3} is given in Section 2 of the supplementary material. It proceeds by establishing $\sqrt{n_I}$-consistency separately for the first-stage and second-stage components of the variance estimator. The proof of Proposition \ref{prop4} follows similar arguments, and is therefore omitted.
Both estimators rely on joint inclusion probabilities, which may be difficult to compute in complex designs. Section \ref{sec:large:entropy} addresses this limitation.

\begin{prop} \label{prop3}
  If assumptions (FS1)-(FS3), (SS1) and (VA1) hold, we have:
    \begin{eqnarray}
      E\left[ A^{-2} n_I \left\{\hat{V}_{HT,A}(\Hat{\tau}_{y})-V_1(\Hat{\tau}_{y})-V_2(\Hat{\tau}_{y})\right\}\right]^2 & = & O(n_I^{-1}), \label{sec:consist:eq5} \\
      E\left[ A^{-2} N_I n_{II0} \left\{\hat{V}_{B}(\Hat{\tau}_{y})-V_3(\Hat{\tau}_{y})\right\}\right]^2 & = & O(n_I^{-1}). \label{sec:consist:eq6}
    \end{eqnarray}
  If in addition assumption (VA2) holds, then:
    \begin{eqnarray}
      \qquad E\left[ A^{-2} n_I \left\{\hat{V}_{HT}(\Hat{\tau}_{y})-V(\Hat{\tau}_{y})\right\}\right]^2 = O(n_I^{-1}) ~~\textrm{and}~~
      \frac{\hat{V}_{HT}(\Hat{\tau}_{y})}{V(\Hat{\tau}_{y})} \rightarrow_{Pr} 1. \label{sec:consist:eq7}
    \end{eqnarray}
\end{prop}

\begin{prop} \label{prop4}
  If assumptions (FS1)-(FS3),(SS1) and (VA1) hold, we have:
    \begin{eqnarray}
      E\left[ A^{-2} n_I \left\{\hat{V}_{YG,A}(\Hat{\tau}_{y})-V_1(\Hat{\tau}_{y})-V_2(\Hat{\tau}_{y})\right\}\right]^2 & = & O(n_I^{-1}). \label{sec:consist:eq8}
    \end{eqnarray}
  If in addition assumption (VA2) holds, then:
    \begin{eqnarray}
      \qquad E\left[ A^{-2} n_I \left\{\hat{V}_{YG}(\Hat{\tau}_{y})-V(\Hat{\tau}_{y})\right\}\right]^2 = O(n_I^{-1}) ~~\textrm{and}~~
      \frac{\hat{V}_{YG}(\Hat{\tau}_{y})}{V(\Hat{\tau}_{y})} \rightarrow_{Pr} 1. \label{sec:consist:eq9}
    \end{eqnarray}
\end{prop}

\subsection{Plug-in estimators}
\label{sec:notation:ssec4}

\noindent We now extend the analysis to smooth scale-invariant functionals of totals. The next result establishes a second-order expansion of the plug-in estimator and shows that the linearization approximation is asymptotically valid. The proof is provided in Appendix \ref{app:proof:prop4b}.

\begin{prop} \label{prop4b}
Suppose that the first-stage assumptions (FS1)-(FS3), the second-stage assumption (SS1), and assumptions (PI1)-(PI2) hold. Then
    \begin{eqnarray}
    E \left(\left[(\hat{\theta}-\theta)-\{g'(\tau_{\bm{\rho}})\}^{\top} \{\hat{\tau}_{\bm{\rho}}-\tau_{\bm{\rho}}\} \right]^2 \right) & = & O(n_I^{-2}), \label{sec:consist:eq10} \\
    E \left[n_I \left| \hat{V}_{HT}(\hat{\theta})-V_p(\hat{\tau}_{l}) \right| \right] & = & O(n_I^{-1/2}), \label{sec:consist:eq11} \\
    E \left[n_I \left| \hat{V}_{YG}(\hat{\theta})-V_p(\hat{\tau}_{l}) \right| \right] & = & O(n_I^{-1/2}), \label{sec:consist:eq12}
    \end{eqnarray}
where $\hat{\tau}_{l}$ is obtained by replacing in \eqref{sec:notation:ssec1:eq8} $\rho(x)$ with the linearized variable $l(x)$.
\end{prop}

\noindent Equation \eqref{sec:consist:eq10} shows that the remainder term in the linearization expansion is of smaller order $O_p(n_I^{-1})$, implying $\sqrt{n_I}$-consistency and asymptotic unbiasedness of the plug-in estimator. Moreover, it guarantees that the linearization variance $V_p(\hat{\tau}_{l})$ is asymptotically equivalent to the true variance. Equations \eqref{sec:consist:eq11} and \eqref{sec:consist:eq12} establish $L^1$-consistency of the linearized variance estimators. This is sufficient for weak consistency and for the validity of normal-based confidence intervals. Stronger mean-square consistency could be obtained under higher moment conditions, at the cost of strengthening (VA1).

\section{High-entropy sampling designs} \label{sec:large:entropy}

\noindent We now show that substantially sharper results can be obtained when the first-stage sampling design has high entropy. In this case, the dependence structure of inclusion indicators is sufficiently regular to (i) eliminate the need for explicit control of joint inclusion probabilities, (ii) weaken moment assumptions, and (iii) establish asymptotic normality of the HT and plug-in estimators. This extends the entropy-based asymptotic theory developed for discrete multistage designs by \citet{chauvet2020inference} to the present continuous-domain setting.

\subsection{Estimators of totals} \label{sec:large:entropy:1}

\noindent We first consider rejective (conditional Poisson) sampling \citep{Haj64}.
Under rejective sampling, higher-order inclusion moments satisfy the bounds in (FS2) automatically \citep{BoiLopRui12,chauvet2020inference}.
Simple random sampling without replacement is obtained as a special case when all first-order inclusion probabilities $\pi_{Ii}$ are equal. Let $p_{rI}(\cdot)$ denote a rejective sampling design on $U_I$ with inclusion probabilities $\pi_{Ii}$, and let $S_{rI}$ be the associated first-stage sample. The H\'ajek variance estimator is defined by
    \begin{align} \label{sec:large:entropy:eq1}
    \Hat{V}_{HAJ}\left(\Hat{\tau}_{\rho}\right) & = \Hat{V}_{HAJ,A}\left(\Hat{\tau}_{\rho}\right) + \Hat{V}_{B}\left(\Hat{\tau}_{\rho}\right)
    \end{align}
where the first-stage component is
    \begin{align} \label{sec:large:entropy:eq2}
    \Hat{V}_{HAJ,A}\left(\Hat{\tau}_{\rho}\right) & = \sum_{u_i \in S_{rI}} (1-\pi_{Ii}) \left(\frac{\Hat{\tau}_{\rho i}}{\pi_{Ii}} - \Hat{\Hat{R}}_{r\pi} \right)^2
    \end{align}
with
    \begin{align*}
    \Hat{\Hat{R}}_{r\pi} = \frac{1}{\hat{d}_{rI}} \sum_{u_i \in S_{rI}} (1-\pi_{Ii}) \frac{\Hat{\tau}_{\rho i}}{\pi_{Ii}}, & \qquad
    \hat{d}_{rI} = \sum_{u_i \in S_{rI}} (1-\pi_{Ii}).
    \end{align*}

\noindent This estimator does not involve the joint inclusion probabilities $\pi_{Iij}$; consequently, assumption (FS3) is no longer required. In contrast to \cite{chauvet2020inference}, truncation of $\hat{d}_{rI}$ is not required here.

\begin{prop} \label{prop5}
  Suppose that rejective sampling is used at the first stage and that assumptions (FS1), (SS1) and (VA1) hold. Then
    \begin{eqnarray} \label{sec:large:entropy:eq3}
      E\left[ A^{-2} n_I \left\{\hat{V}_{HAJ,A}(\Hat{\tau}_{\rho})-V_1(\Hat{\tau}_{\rho})-V_2(\Hat{\tau}_{\rho})\right\}\right]^2 & = & o(1).
    \end{eqnarray}
  If, in addition, (VA2) holds, then
    \begin{eqnarray} \label{sec:large:entropy:eq4}
      \frac{\Hat{\tau}_{\rho}-\tau_{\rho}}{\sqrt{V(\Hat{\tau}_{\rho})}} & \longrightarrow_{\mathcal{L}} & \mathcal{N}(0,1),
    \end{eqnarray}
  where $\rightarrow_{\mathcal{L}}$ stands for the convergence in distribution. Also:
    \begin{align} \label{sec:large:entropy:eq5}
      E\left[ A^{-2} n_I \left\{\hat{V}_{HAJ}(\Hat{\tau}_{\rho})-V(\Hat{\tau}_{\rho})\right\}\right]^2 = o(1),
      & \qquad \frac{\hat{V}_{HAJ}(\Hat{\tau}_{\rho})}{V(\Hat{\tau}_{\rho})} \rightarrow_{Pr} 1.
    \end{align}
\end{prop}

\noindent Proposition \ref{prop5} shows that, under rejective sampling, the standardized HT estimator converges to a Gaussian limit without requiring explicit control of second-order inclusion probabilities. This result parallels classical finite-population central limit theorems, but in the present continuous-domain setting the proof requires careful control of the second-stage variance components. The proof is given in Appendix \ref{app:proof:prop5}. \\
\noindent The previous result is not restricted to rejective sampling. It extends to any fixed-size design that is asymptotically close to rejective sampling in chi-square distance. For a design $p_I(\cdot)$ with inclusion probabilities $\pi_{Ii}$, define
    \begin{eqnarray} \label{sec:large:entropy:eq6}
      d_2(p_I,p_{rI})=\sum_{s_I \subset U_I;~p_{rI}(s_I)>0} \frac{\left\{p_I(s_I)-p_{rI}(s_I)\right\}^2}{p_{rI}(s_I)}.
    \end{eqnarray}
The design $p_I(\cdot)$ is said to be close to $p_{rI}(\cdot)$ if $d_2(p_I,p_{rI}) \to 0$; this condition is satisfied, for example, by the Rao-Sampford design \citep{Sam67}. Proposition \ref{prop6} extends Proposition \ref{prop5} to high-entropy sampling designs; see Appendix \ref{app:proof:prop6} for a proof.

\begin{prop} \label{prop6}
  Suppose that (FS1), (SS1), (VA1)-(VA2) hold, and that $d_2(p_I,p_{rI}) \to 0$. Then
    \begin{eqnarray} \label{sec:large:entropy:eq7}
      \frac{\hat{\tau}_{\rho}-\tau_{\rho}}{\sqrt{V(\hat{\tau}_{\rho})}} & \longrightarrow_{\mathcal{L}} & \mathcal{N}(0,1).
    \end{eqnarray}
  Also:
    \begin{align} \label{sec:large:entropy:eq8}
      E\left[ A^{-2} n_I \left| \hat{V}_{HAJ}(\Hat{\tau}_{\rho})-V(\Hat{\tau}_{\rho})\right| \right] = o(1),
      & \qquad \frac{\hat{V}_{HAJ}(\Hat{\tau}_{\rho})}{V(\Hat{\tau}_{\rho})} \rightarrow_{Pr} 1.
    \end{align}
\end{prop}

\subsection{Plug-in estimators} \label{sec:large:entropy:2}

\noindent We now extend the asymptotic theory for high-entropy designs to smooth functionals of totals. Proposition \ref{prop4b} established that, under general conditions, the plug-in estimator admits a valid linearization expansion.
We now show that, under high-entropy first-stage designs, this expansion leads to asymptotic normality and consistent variance estimation without requiring joint inclusion probabilities.

\begin{prop} \label{prop7}
  Suppose that (FS1), (SS1), (PI1)-(PI2) hold, and that $d_2(p_I,p_{rI}) \to 0$. Suppose that the first-stage variance component of $\hat{\theta}$ does not degenerate, i.e. there exists some constant $m_3$ such that
         \begin{align} \label{prop7:eq1}
         n_I V_1(\Hat{\theta}) \geq m_3 >0.
         \end{align}
  Then
    \begin{eqnarray} \label{prop7:eq2}
      \frac{\hat{\theta}-\theta}{\sqrt{V(\hat{\theta})}} & \longrightarrow_{\mathcal{L}} & \mathcal{N}(0,1).
    \end{eqnarray}
  Also:
    \begin{align} \label{prop7:eq3}
      E\left[ A^{-2} n_I \left| \hat{V}_{HAJ}(\Hat{\theta})-V(\Hat{\theta})\right| \right] = o(1),
      & \qquad \frac{\hat{V}_{HAJ}(\Hat{\theta})}{V(\Hat{\theta})} \rightarrow_{Pr} 1.
    \end{align}
\end{prop}

\noindent Proposition \ref{prop7} shows that, under high-entropy first-stage designs, the plug-in estimator inherits the asymptotic Gaussian behavior of its linearized counterpart. The additional non-degeneracy condition \eqref{prop7:eq1} ensures that the first-stage variability remains asymptotically dominant. Combined with the $L^1$-consistency of the H\'ajek variance estimator, this result justifies the use of normal-approximation confidence intervals for smooth functionals of totals in the continuous-domain two-stage framework. The proof is provided in Appendix \ref{app:proof:prop7}.

\section{Simulation study} \label{sec:simulation}

\subsection{Simulation set-up}

\begin{figure}[!hbtp]
    \centering
    \includegraphics[width=0.5\linewidth]{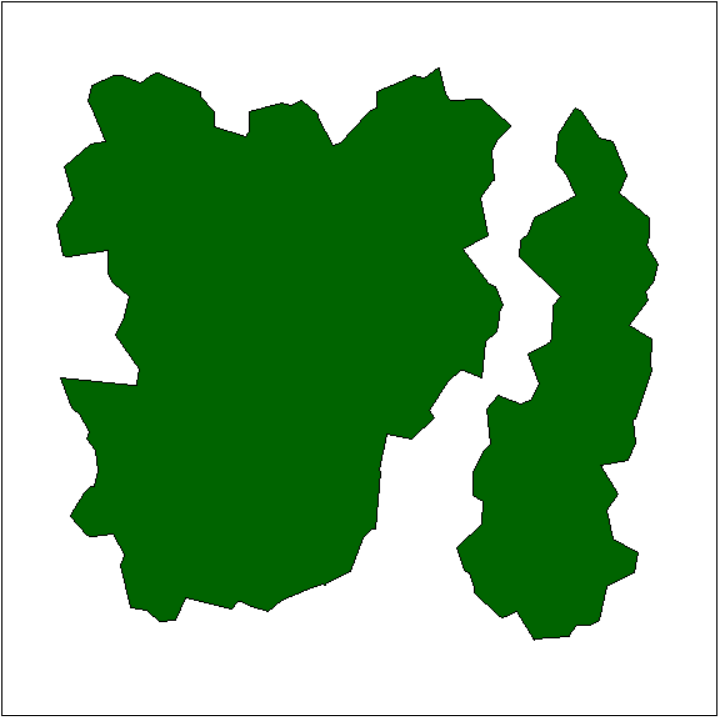}
    \caption{The study territory including the forest area (in green)}
    \label{fig:territory}
\end{figure}

\noindent To evaluate the sampling designs and estimators considered in this article, we conduct a simulation study in a realistic  artificial forest population. The study territory $\mcF$ is a square of area $A=100$ km$^{2}$. The forest occupies part of this territory and contains a total of $N=228,385$ trees, with an overall wood volume of $V = 158,534.74$ m$^3$. For each tree, the available attributes are species, height, breast-height diameter, and volume. The spatial configuration of the forest within $\mcF$ is displayed in Figure \ref{fig:territory}. The dataset is available at \citet{duong23art}, version 2.

\begin{figure}[!hbtp]
    \centering
    \footnotesize
    \begin{tabular}{cc}
        \includegraphics[width=0.5\linewidth]{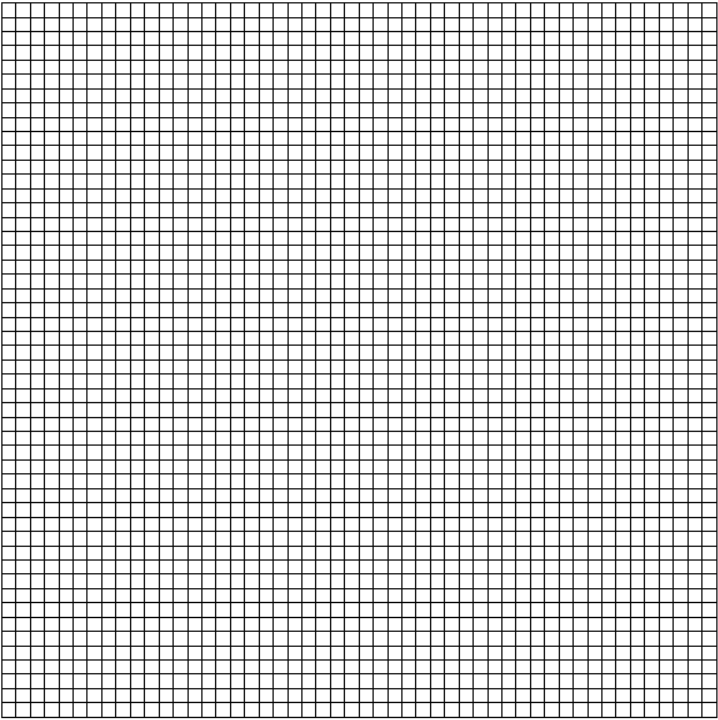} & \includegraphics[width=0.5\linewidth]{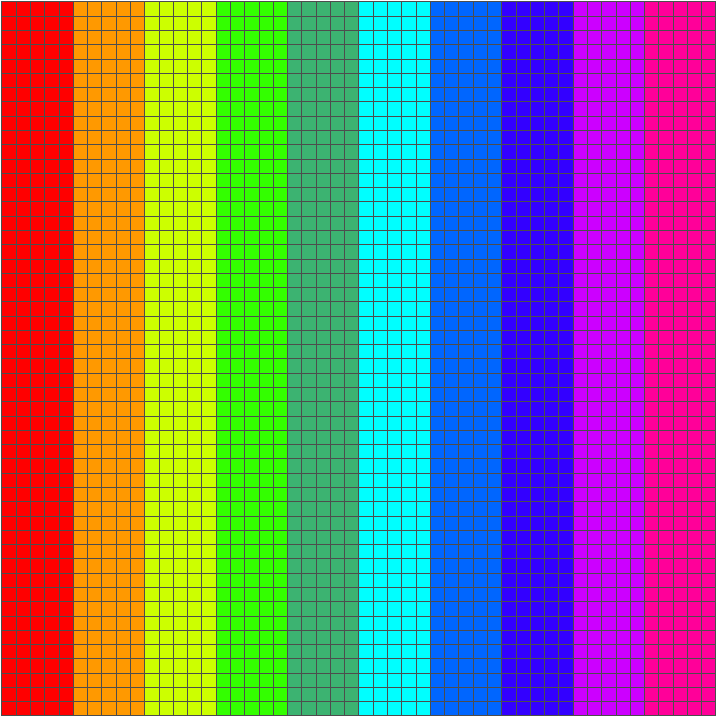}\\
        a & b \\
        \multicolumn{2}{c}{\includegraphics[width=1\linewidth]{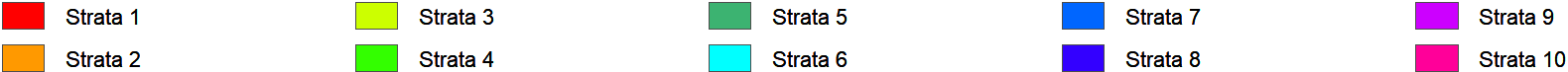}}
    \end{tabular}
    \caption{The territory with a regular tesselation (a) and partitioned into 10 strata (b)}
    \label{fig:regular_grid}
\end{figure}

\begin{figure}[!hbtp]
    \centering
    \footnotesize
    \begin{tabular}{cc}
    \includegraphics[width=0.5\linewidth]{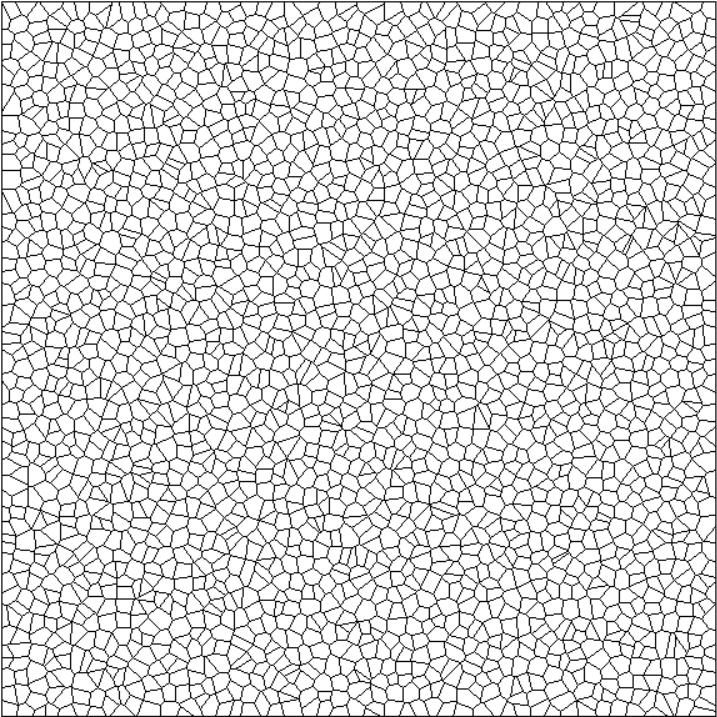} & \includegraphics[width=0.5\linewidth]{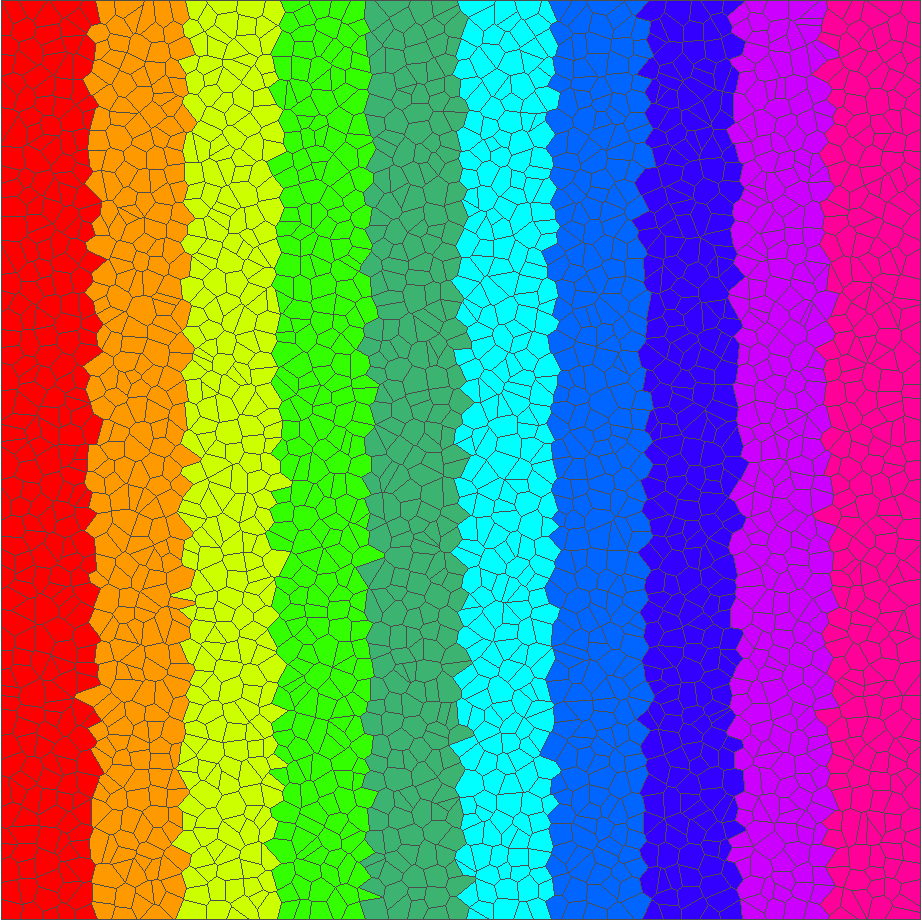} \\
        a & b \\
        \multicolumn{2}{c}{\includegraphics[width=1\linewidth]{Images/legend_stratification__4.png}}
    \end{tabular}
    \caption{The territory with a irregular tesselation (a) and partitioned into 10 strata (b)}
    \label{fig:irregular_grid}
\end{figure}

\noindent The territory $\mcF$ is partitioned into PSUs to form the first-stage population $U_I$. Two configurations are considered. Firstly, a regular grid with $N_I=2,500$ square cells of equal area $A_i=A_0=40,000$ m$^{2}$. Secondly, an irregular partition into $N_I=2,500$ cells with heterogeneous areas ranging from approximately $11,954$ m$^2$ to $78,185$ m$^2$. These configurations are illustrated in Figures \ref{fig:regular_grid}a and \ref{fig:irregular_grid}a, respectively.

\noindent For the first-stage sample $S_I$ of size $n_I$, we consider: simple random sampling without replacement (SRSwor), stratified simple random sampling with proportional allocation (StrSRS), and rejective sampling (REJ) with probabilities proportional to PSU area. For StrSRS, the PSUs are grouped into 10 strata (Figures \ref{fig:regular_grid}b and \ref{fig:irregular_grid}b). Note that under the regular grid, all PSUs have equal area, so rejective sampling reduces to simple random sampling. We consider first-stage sample sizes $n_I \in \{100,200,500\}$. At the second stage, for each selected PSU $u_i \in S_I$, we independently select $n_{IIi}=2$ points independently by uniform sampling over the cell. \\
\noindent For each selected point, a circular plot of area $A_r=500$ m$^2$ is centered at the point, and all trees whose centers lie inside the plot are measured. The tree-level measurements $y_k$ are converted into a local attribute density $\rho(x)$ at the sampled point, as defined in equation \eqref{weight:share}; see also \cite{chaboubri23,bouriaud2024wsm} for details. We consider estimation of integral totals of the form $\tau_{\rho}$, corresponding to total volume of wood $V$, and on the other hand to the total number of trees $N$. These are estimated using the HT estimator $\hat{\tau}_{\rho}$ defined in equation \eqref{sec:notation:ssec1:eq8}. Variance estimation is performed using $\hat{V}_{HT}(\hat{\tau}_{\rho})$, $\hat{V}_{YG}(\hat{\tau}_{\rho})$, and $\hat{V}_{HAJ}(\hat{\tau}_{\rho})$, defined in equations \eqref{sec:notation:ssec1:eq12}, \eqref{sec:notation:ssec1:eq14}, and \eqref{sec:large:entropy:eq1}, respectively. \\
\noindent We also consider ratios of the form $\theta = \tau_{\rho_1}/\tau_{\rho_2}$, corresponding to percentage of tree volume per species ($R_v$), and to percentage of total number of trees per species ($R_N$). The dataset contains $33$ species, including conifers (e.g., \textit{Abies alba}, \textit{Picea abies}, \textit{Pinus sylvestris}) and broadleaved species (e.g., \textit{Fagus sylvatica}, \textit{Castanea sativa}, \textit{Fraxinus excelsior}). We focus on \textit{Abies alba}, which comprises $N_{abi}=83,395$ trees and a total volume of $V_{abi} = 72,531.34$ m$^3$, corresponding to $R_v = 46\%$ and $R_N = 37\%$. Let $\Hat{R}_v$ and $\Hat{R}_N$ denote the corresponding ratio estimators, using the plug-in estimator $\hat{R} = \hat{\theta} = \hat{\tau}_{\rho_1}/\hat{\tau}_{\rho_2}$ as defined in section \ref{sec:notation:ssec2}. The associated linearized variable is
    \begin{equation} \label{lin:ratio}
    \hat{l}(x) = \frac{1}{\hat{\tau}_{\rho_2}} \left\{\rho_1(x)-\hat{\theta} \rho_2(x)\right\},
    \end{equation}
which is substituted into the variance estimators. \\
\noindent For each population configuration, the sampling and estimation procedure is repeated $M=50,000$ times. For an estimator $\hat{\theta}$, we compute the Monte Carlo percent relative bias
    \begin{equation} \label{eq:RB_MC}
    {RB}_{MC}(\hat{\theta}) = \frac{M^{-1} \sum_{r=1}^M \left(\hat{\theta}^m - \theta\right)}{\theta} \times 100,
    \end{equation}
with $\hat{\theta}^m$ the value of the estimator on the $m-th$ sample.
For variance estimators, the true variance is approximated by an independent Monte Carlo run of $50,000$ replications. Also, we compute the Monte Carlo percent relative stability
    \begin{equation}
    {RS}_{MC}(\hat{V}) = \frac{\left[M^{-1} \sum_{r=1}^M \left\{\hat{V}^r - V(\hat{\theta}) \right\}^2\right]^{0.5}}{V(\hat{\theta})} \times 100.
    \end{equation}
Finally, we compute empirical confidence interval error rates based on asymptotic normality, using a nominal one-tailed error rate of $2.5 \%$ in each tail (corresponding to a $95\%$ confidence level). These are denoted by $CI_{MC}(\hat{V})$.

\subsection{Simulation results}

\noindent The simulation results are presented in Tables \ref{tab:total_est}–\ref{tab:ratio_esp_trees}. For SRSwor and StrSRS, the HT and YG variance estimators coincide exactly. In the case of rejective sampling (REJ), their numerical results were virtually identical. For the sake of brevity, we therefore report only the results obtained with the HT variance estimator. \\
\noindent The results for the estimation of total wood volume are displayed in Table \ref{tab:total_est}. As expected, the HT total estimator is unbiased across all sampling designs and sample sizes. The variance decreases steadily as the first-stage sample size $n_I$ increases, in agreement with the theoretical order of magnitude established earlier. The HT variance estimator is essentially unbiased in all scenarios. The HAJ variance estimator exhibits a slight negative bias for small sample sizes under StrSRS; however, this bias diminishes as $n_I$ increases. In terms of relative stability, the HT and HAJ variance estimators perform similarly. \\
\noindent The empirical coverage rates of the normality-based confidence intervals are close to the nominal 95 \% level, which supports the asymptotic results derived in the previous sections. A slightly lower coverage rate is observed for StrSRS combined with the HAJ variance estimator, which is consistent with the small negative bias noted above.

\begin{table}[!hbtp]
\footnotesize
\centering
\begin{tabular}{lrrrrrrrrr}
\toprule
Design & Sample & $RB_{MC}$ & $var_{MC}$ & \multicolumn{2}{c}{$RB_{MC}$} & \multicolumn{2}{c}{$RS_{MC}$} & \multicolumn{2}{c}{$CI_{MC}$} \\
& size & $(\Hat{\tau}_{v})$ & $(\Hat{\tau}_{v})$
& $\Hat{V}_{HT}$ & $\Hat{V}_{HAJ}$
& $\Hat{V}_{HT}$ & $\Hat{V}_{HAJ}$
& $\Hat{V}_{HT}$ & $\Hat{V}_{HAJ}$ \\
& & (\%) & ($\times 10^8$) & (\%) & (\%) & (\%) & (\%) & (\%) & (\%) \\
\midrule
\multicolumn{10}{c}{\textit{Regular grid}} \\
\addlinespace
SRSwor & 100 & -0.02 & 9.59 & 0.51 & -0.48 & 33.73 & 33.40 & 0.93 & 0.93 \\
       & 200 & -0.19 & 4.66 & 0.63 & 0.14 & 23.65 & 23.53 & 0.94 & 0.94 \\
       & 500 & 0.00 & 1.72 & 1.25 & 1.06 & 14.59 & 14.55 & 0.95 & 0.95 \\
\addlinespace
StrSRS & 100 & -0.03 & 8.81 & -0.38 & -10.16 & 33.76 & 32.08 & 0.93 & 0.92 \\
       & 200 & -0.06 & 4.30 & -0.38 & -5.17 & 23.54 & 22.98 & 0.94 & 0.93 \\
       & 500 & -0.03 & 1.60 & -0.15 & -1.95 & 14.33 & 14.20 & 0.95 & 0.94 \\
\midrule
\multicolumn{10}{c}{\textit{Irregular grid}} \\
\addlinespace
SRSwor & 100 & -0.05 & 10.22 & -0.43 & -1.41 & 36.81 & 36.47 & 0.93 & 0.93 \\
       & 200 & -0.08 & 5.03 & -1.30 & -1.78 & 25.66 & 25.56 & 0.94 & 0.94 \\
       & 500 & -0.00 & 1.84 & 0.08 & -0.10 & 15.77 & 15.74 & 0.95 & 0.95 \\
\addlinespace
StrSRS & 100 & -0.06 & 9.28 & 0.29 & -9.55 & 37.08 & 34.75 & 0.93 & 0.92 \\
       & 200 & -0.06 & 4.52 & 0.75 & -4.10 & 26.09 & 25.15 & 0.94 & 0.93 \\
       & 500 & 0.02 & 1.73 & -1.69 & -3.46 & 15.68 & 15.69 & 0.94 & 0.94 \\
\addlinespace
REJ    & 100 & 0.02 & 9.67 & -0.60 & -1.57 & 33.50 & 33.18 & 0.93 & 0.93 \\
       & 200 & 0.10 & 4.69 & -0.10 & -0.58 & 23.42 & 23.31 & 0.94 & 0.94 \\
       & 500 & 0.02 & 1.72 & 0.31 & 0.13 & 14.38 & 14.35 & 0.95 & 0.95 \\
\bottomrule
\end{tabular}
\caption{Monte Carlo results for the estimation of total wood volume $\hat{\tau}_v$: relative bias, variance, and performance of variance estimators (RB, RS, CI).}
\label{tab:total_est}
\end{table}

\noindent The results for estimating the total number of trees are presented in Table \ref{tab:number_trees}. The qualitative conclusions are the same as for total wood volume: the estimator is unbiased, the variance decreases with increasing sample size, and the variance estimators perform comparably, with coverage rates close to the nominal level.

\begin{table}[!hbtp]
    \centering
    \footnotesize
        \begin{tabular}{lrrrrrrrrr}
        \toprule
        Design & Sample & $RB_{MC}$ & $var_{MC}$ & \multicolumn{2}{c}{$RB_{MC}$} & \multicolumn{2}{c}{$RS_{MC}$} & \multicolumn{2}{c}{$CI_{MC}$} \\
        & size & $(\Hat{N})$ & $(\Hat{N})$ & $\Hat{V}_{HT} $ & $\Hat{V}_{HAJ}$ &  $\Hat{V}_{HT}$ & $\Hat{V}_{HAJ}$ & $\Hat{V}_{HT}$ & $\Hat{V}_{HAJ}$ \\
         & & (\%) & ($\times$ 10$^8$) & (\%) & (\%) & (\%) & (\%) & (\%) & (\%) \\
         \midrule
         \multicolumn{10}{c}{\textit{Regular grid}} \\
         \addlinespace
         SRSwor & 100 & -0.03 & 14.69 & -0.06 & -1.05 & 24.42 & 24.20 & 0.94 & 0.94\\
         & 200 & -0.14 & 7.15 & -0.63 & -1.21 & 17.07 & 17.01 & 0.94 & 0.94  \\
         & 500 & -0.01 & 2.54 & 0.73 & 0.54 & 14.59 & 14.55 & 0.95 & 0.95  \\
         \addlinespace
         StrSRS & 100 & -0.11 & 13.16 & 0.52 & -9.44 & 25.27 & 24.63 & 0.94 & 0.92  \\
         & 200 & -0.06 & 6.36 & 0.79 & -4.15 & 17.31 & 16.95 & 0.95 & 0.94 \\
         & 500 & -0.01 & 2.32 & -0.14 & -2.03 & 10.23 & 10.23 & 0.95 & 0.94 \\
         \midrule
         \multicolumn{10}{c}{\textit{Irregular grid}} \\
         \addlinespace
         SRSwor & 100 & 0.03 & 15.78 & -0.54 & -1.52 & 28.59 & 28.34 & 0.94 & 0.93 \\
         & 200 & -0.05 & 7.62 & -0.44 & -0.93 & 19.93 & 19.85 & 0.94 & 0.94 \\
         & 500 & 0.01 & 2.72 & 0.30 & 0.11 & 11.93 & 11.91 & 0.95 & 0.95 \\
         \addlinespace
         StrSRS & 100 & 0.04 & 14.35 & -0.67 & -10.51 & 29.27 & 28.37 & 0.94 & 0.92  \\
         & 200 & -0.06 & 6.84 & 0.77 & -4.17 & 20.38 & 19.80 & 0.94 & 0.94 \\
         & 500 & 0.01 & 2.52 & -1.27 & -3.13 & 12.05 & 12.16 & 0.95 & 0.94 \\
         \addlinespace
         REJ & 100 & -0.02 & 14.66 & -0.22 & -1.21 & 24.42 & 24.21 & 0.94 & 0.94 \\
         & 200 & 0.12 & 7.14 & -0.79 & -1.28 & 16.97 & 16.91 & 0.94 & 0.94 \\
         & 500 & 0.01 & 2.53 & -0.14 & -0.32 & 10.21 & 10.19 & 0.95 & 0.95 \\
         \bottomrule
    \end{tabular}
    \caption{Monte Carlo results for the estimation of total number of trees $\hat{N}$: relative bias, variance, and performance of variance estimators (RB, RS, CI).}
    \label{tab:number_trees}
\end{table}

\noindent The results for the estimation of the ratio of total wood volume attributable to \textit{Abies alba} are reported in Table \ref{tab:ratio_esp_volume}. The estimator of the ratio is essentially unbiased across all sampling designs and sample sizes. Its variance decreases as the first-stage sample size $n_I$ increases, in agreement with the theoretical results established in Proposition \ref{prop4b}. The HT linearized variance estimator shows a small negative bias for small sample sizes, but this bias diminishes as $n_I$ increases, which is consistent with the asymptotic theory. The HAJ variance estimator, however, exhibits a substantial bias under StrSRS. In terms of relative stability, the two variance estimators perform similarly.

\begin{table}[!hbtp]
\centering
\footnotesize
    \begin{tabular}{lrrrrrrrrr}
        \toprule
        Design & Sample & $RB_{MC}$ & $var_{MC}$ & \multicolumn{2}{c}{$RB_{MC}$} & \multicolumn{2}{c}{$RS_{MC}$} & \multicolumn{2}{c}{$CI_{MC}$} \\
        & size & $(\Hat{R}_{v})$ & $(\Hat{R}_{v})$ & $\Hat{V}_{HT}$ & $\Hat{V}_{HAJ}$ &  $\Hat{V}_{HT}$ & $\Hat{V}_{HAJ}$ & $\Hat{V}_{HT}$ & $\Hat{V}_{HAJ}$ \\
         & & (\%) & ($\times$ 10$^{-3}$) & (\%) & (\%) & (\%) & (\%) & (\%) & (\%) \\
         \midrule
         \multicolumn{10}{c}{\textit{Regular grid}} \\
         \addlinespace
         SRSwor & 100 & -0.25 & 5.10 & -7.67 & -8.56 & 31.76 & 31.70 & 0.92 & 0.92\\
         & 200 & -0.13 & 2.51 & -3.13 & -3.58 & 22.58 & 22.54 & 0.94 & 0.94  \\
         & 500 & -0.06 & 0.98 & -1.75 & -1.91 & 13.38 & 13.38 & 0.94 & 0.94  \\
         \addlinespace
         StrSRS & 100 & -0.32 & 5.04 & -6.47 & -15.51 & 31.43 & 31.76 & 0.92 & 0.91  \\
         & 200 & -0.19 & 2.52 & -4.00 & -8.49 & 22.41 & 22.63 & 0.93 & 0.93 \\
         & 500 & -0.07 & 0.99 & -2.55 & -4.17 & 13.25 & 13.43 & 0.94 & 0.94 \\
         \midrule
         \multicolumn{10}{c}{\textit{Irregular grid}} \\
         \addlinespace
         SRSwor & 100 & -0.23 & 5.33 & -8.61 & -9.49 & 33.07 & 33.00 & 0.92 & 0.91 \\
         & 200 & -0.21 & 2.62 & -3.74 & -4.19 & 23.72 & 23.68 & 0.93 & 0.93 \\
         & 500 & -0.03 & 1.03 & -2.19 & -2.36 & 14.18 & 14.19 & 0.94 & 0.94 \\
         \addlinespace
         StrSRS & 100 & -0.33 & 5.22 & -7.02 & -16.01 & 32.85 & 33.05 & 0.92 & 0.91 \\
         & 200 & -0.11 & 2.65 & -5.31 & -9.72 & 23.37 & 23.73 & 0.93 & 0.93 \\
         & 500 & -0.05 & 1.02 & -1.92 & -3.55 & 14.05 & 14.12 & 0.94 & 0.94 \\
         \addlinespace
         REJ & 100 & -0.27 & 5.10 & -7.68 & -8.57 & 31.76 & 31.69 & 0.92 & 0.92 \\
         & 200 & -0.10 & 2.50 & -3.29 & -3.74 & 22.54 & 22.50 & 0.93 & 0.93 \\
         & 500 & -0.06 & 0.98 & -1.91 & -2.07 & 13.50 & 13.50 & 0.94 & 0.94 \\
         \bottomrule
    \end{tabular}
\caption{Monte Carlo results for the estimation of the proportion of wood volume of \textit{Abies alba} $\hat{R}_{v}$: relative bias, variance, and performance of variance estimators (RB, RS, CI).}
\label{tab:ratio_esp_volume}
\end{table}

\noindent The empirical coverage rates are generally close to the nominal level. As expected, slightly lower coverage is observed in settings where the variance estimator is negatively biased. The simulation results for the estimation of the percentage of \textit{Abies alba} in terms of number of trees are presented in Table \ref{tab:ratio_esp_trees}. The qualitative conclusions are the same as for the volume ratio.

\begin{table}[!hbtp]
\centering
\footnotesize
     \begin{tabular}{lrrrrrrrrr}
        \hline
        Design & Sample & $RB_{MC}$ & $var_{MC}$ & \multicolumn{2}{c}{$RB_{MC}$} & \multicolumn{2}{c}{$RS_{MC}$} & \multicolumn{2}{c}{$CI_{MC}$} \\
        & size & $(\Hat{R}_{N})$ & $(\Hat{R}_{N})$ & $\Hat{V}_{HT}$ & $\Hat{V}_{HAJ}$ &  $\Hat{V}_{HT}$ & $\Hat{V}_{HAJ}$ & $\Hat{V}_{HT}$ & $\Hat{V}_{HAJ}$ \\
         & & (\%) & ($\times$ 10$^{-3}$) & (\%) & (\%) & (\%) & (\%) & (\%) & (\%) \\
         \toprule
         \multicolumn{10}{c}{\textit{Regular grid}} \\
         \addlinespace
         SRSwor & 100 & 0.06 & 2.33 & -5.35 & -6.27 & 32.86 & 32.70 & 0.93 & 0.93 \\
         & 200 & 0.04 & 1.14 & -1.72 & -2.18 & 23.03 & 22.96 & 0.94 & 0.94  \\
         & 500 & -0.02 & 0.45 & -1.21 & -1.37 & 13.24 & 13.24 & 0.95 & 0.95  \\
         \addlinespace
         StrSRS & 100 & -0.02 & 2.31 & -4.87 & -14.06 & 32.60 & 32.26 & 0.93 & 0.91  \\
         & 200 & -0.04 & 1.13 & -1.81 & -6.39 & 22.76 & 22.50 & 0.94 & 0.93 \\
         & 500 & -0.05 & 0.45 & -1.43 & -3.07 & 13.07 & 13110 & 0.94 & 0.94 \\
         \midrule
         \multicolumn{10}{c}{\textit{Irregular grid}} \\
         \addlinespace
         SRSwor & 100 & 0.17 & 2.43 & -4.91 & -5.83 & 34.74 & 34.55 & 0.93 & 0.93 \\
         & 200 & 0.02 & 1.21 & -2.96 & -3.41 & 24.64 & 24.58 & 0.94 & 0.94 \\
         & 500 & 0.03 & 0.47 & -1.04 & -1.21 & 14.25 & 14.24 & 0.95 & 0.95 \\
         \addlinespace
         StrSRS & 100 & -0.01 & 2.44 & -5.87 & -14.96 & 34.46 & 34.04 & 0.93 & 0.91 \\
         & 200 & 0.01 & 1.20 & -2.95 & -7.48 & 24.15 & 23.98 & 0.94 & 0.93 \\
         & 500 & 0.00 & 0.46 & 0.07 & -1.59 & 14.26 & 14.08 & 0.95 & 0.94 \\
         \addlinespace
         REJ & 100 & 0.08 & 2.31 & -4.01 & -4.94 & 33.15 & 32.95 & 0.93 & 0.93 \\
         & 200 & 0.02 & 1.14 & -1.84 & -2.29 & 22.82 & 22.75 & 0.94 & 0.94 \\
         & 500 & 0.01 & 0.45 & -1.54 & -1.71 & 13.20 & 13.20 & 0.94 & 0.94 \\
         \bottomrule
    \end{tabular}
\caption{Monte Carlo results for the estimation of the proportion of trees of \textit{Abies alba} $\hat{R}_{N}$: relative bias, variance, and performance of variance estimators (RB, RS, CI).}
\label{tab:ratio_esp_trees}
\end{table}

\section{Conclusion} \label{sec:conclusion}

\noindent This paper develops a design-based asymptotic theory for two-stage sampling over a continuous spatial domain. Within a fixed-area framework, we establish consistency and asymptotic normality of Horvitz–Thompson estimators, together with the consistency of associated variance estimators. The results are further extended to smooth functions of totals through linearization arguments. Our analysis complements existing asymptotic theory for multistage sampling in finite populations by explicitly accounting for the continuous nature of second-stage units. In particular, we provide conditions under which standard variance estimators remain valid. \\
\noindent These results are motivated by applications in environmental surveys, such as national forest inventories, where two-stage sampling over continuous territories is routinely implemented. In such settings, extensions to more complex designs, including two-phase sampling with post-stratification, are of practical importance. While convergence results for two-phase designs have been established in the environmental literature \citep{fattorini2017design}, the extension to nonlinear functionals and the derivation of asymptotic normality for continuous populations remain open problems. Addressing these questions would provide a unified asymptotic framework for modern survey designs combining multistage, multiphase, and spatial sampling features.

\section*{Supplementary Material}

\noindent The supplementary material contains additional proofs of intermediary results.

\section*{Acknowledgments}

This work was supported by a grant from the institute of Mathematics for Planet Earth (iMPT, project CONIFER), and by the France 2030 program "Initiative d'Excellence Lorraine (LUE)", reference ANR-15-IDEX-04-LUE, project Artemis-BOOSTFOR. This work was also supported by project "Interdisciplinary Cloud and Big Data Center at Stefan cel Mare University of Suceava" [POC/398/1/1, 343/390019], co-founded by the European Union. OB acknowledges funding from EU Horizon Europe Research and Innovation Program, project INFORMA Grant agreement No. 101060309. \\
ChatGPT was used in this paper as a language assistant for improving the writing of this paper and for limited code editing.

\appendix


\section{Additional notations} \label{app:add:not}

\subsection{Variance decomposition} \label{app:add:not:1}

The variance components of $\Hat{\tau}_{\rho}$ are defined as
    \begin{align} 
    V(\Hat{\tau}_{\rho}) & = \sum_{u_i \in U_I} \sum_{u_j \in U_I} \Delta_{Iij} \frac{\tau_{\rho i}}{\pi_{Ii}} \frac{\tau_{\rho j}}{\pi_{Ij}} + \sum_{u_i \in U_I} \frac{1-\pi_{Ii}}{\pi_{Ii}} V_i + \sum_{u_i \in U_I} V_i, \nonumber \\
    & = V_1(\Hat{\tau}_{\rho})+V_2(\Hat{\tau}_{\rho})+V_3(\Hat{\tau}_{\rho}),
    \end{align}
where $\Delta_{Iij}=\pi_{Iij}-\pi_{Ii}\pi_{Ij}$, and
    \begin{equation} \label{sec:notation:ssec1:eq10}
    V_i = \int_{x \in u_i} \frac{\left\{\rho(x)\right\}^2}{\pi_{IIi}(x)} \D{x}
    + \int_{x \in u_i} \int_{x' \in u_i} \left\{\pi_{IIi}(x,x')-\pi_{IIi}(x)\pi_{IIi}(x')\right\} \frac{\rho(x)}{\pi_{IIi}(x)}
    \frac{\rho(x')}{\pi_{IIi}(x')} \D{x} \D{x'},
    \end{equation}
see Proposition 1 in \citet{duong24}. Since the $n_{IIi}$ points are independently selected within $u_i$ according to the marginal density $f_{IIi}(\cdot)$, this term simplifies to
    \begin{align} \label{sec:notation:ssec1:eq11}
    V_i & = \int_{x \in u_i} \pi_{IIi}(x) \left\{\frac{\rho(x)}{\pi_{IIi}(x)}-\frac{\tau_{\rho i}}{n_{IIi}}\right\}^2 \D{x}.
    \end{align}

\subsection{Variance estimation} \label{app:add:not:2}

\noindent The Horvitz-Thompson variance estimator is defined by
    \begin{align} 
    \Hat{V}_{HT}\left(\Hat{\tau}_{\rho}\right) & = \sum_{u_i \in S_I} \sum_{u_j \in S_I} \frac{\Delta_{Iij}}{\pi_{Iij}} \frac{\Hat{\tau}_{\rho i}}{\pi_{Ii}} \frac{\Hat{\tau}_{\rho j}}{\pi_{Ij}} + \sum_{u_i \in S_I} \frac{\Hat{V}_i}{\pi_{Ii}} \nonumber \\
    & = \Hat{V}_{HT,A}\left(\Hat{\tau}_{\rho}\right) + \Hat{V}_{B}\left(\Hat{\tau}_{\rho}\right)
    \end{align}
where
    \begin{align} \label{sec:notation:ssec1:eq13}
    \Hat{V}_i & = \sum_{u_i \in S_I} \left\{\frac{\rho(x)}{\pi_{IIi}(x)}\right\}^2
    + \sum_{x \neq x' \in S_{I}}
    \frac{\pi_{IIi}(x,x')-\pi_{Ii}(x) \pi_{Ii}(x')}{\pi_{IIi}(x,x')}
    \frac{\rho(x)}{\pi_{Ii}(x)} \frac{\rho(x')}{\pi_{Ii}(x')} \nonumber \\
    & = \frac{n_{IIi}}{n_{IIi}-1} \sum_{x \in S_{IIi}} \left\{\frac{\rho(x)}{\pi_{IIi}(x)}-\frac{\Hat{\tau}_{\rho i}}{n_{IIi}}\right\}^2.
    \end{align}

\noindent The estimator $\Hat{V}_{HT}\left(\Hat{\tau}_{\rho}\right)$ is design-unbiased for $V\left(\Hat{\tau}_{\rho}\right)$, provided that $\pi_{Iij}>0$ for any $u_i,u_j \in U_I$ and that $n_{IIi} \geq 2$ for any $u_i \in U_I$; see Proposition 2 in \citet{duong24}. \\

\noindent The Yates-Grundy variance estimator is defined by
    \begin{align} 
    \Hat{V}_{YG}\left(\Hat{\tau}_{\rho}\right) & = -\frac{1}{2}\sum_{u_i \neq u_j \in S_I} \frac{\Delta_{Iij}}{\pi_{Iij}} \left(\frac{\Hat{\tau}_{\rho i}}{\pi_{Ii}}- \frac{\Hat{\tau}_{\rho j}}{\pi_{Ij}}\right)^2 + \sum_{u_i \in S_I} \frac{\Hat{V}_i}{\pi_{Ii}} \nonumber \\
    & = \Hat{V}_{YG,A}\left(\Hat{\tau}_{\rho}\right) + \Hat{V}_{B}\left(\Hat{\tau}_{\rho}\right).
    \end{align}

\section{Preliminary results} \label{app:prelim}

\noindent The proofs of lemmas are given in Section 1 of the supplementary material.

\subsection{Inequalities related to the first-stage inclusion probabilities} \label{app:prelim:sec1}

\begin{lemma} \label{lem0}
Suppose that Assumption (FS1) holds, and that the sample $S_I$ is of fixed-size $n_I$. Then
    \begin{align*}
    d_I = \sum_{u_i \in U_I} \pi_{Ii}(1-\pi_{Ii}) & \geq c_{I1} (1-f_{I0}) n_I, \\
    \hat{d}_I = \sum_{u_i \in S_I} (1-\pi_{Ii}) & \geq \frac{1}{2} c_{I1} (1-f_{I0}) n_I. 
    \end{align*}
\end{lemma}

\begin{lemma} \label{lem0b}
Suppose that Assumptions (FS1) and (FS2) hold. Then some constants $C'_{I3},C'_{I4},C'_{I5}$ exist such that for any $i \neq j \neq i' \neq j' \in U_I$:
    \begin{align}
    \pi_{Iij} & \leq \left\{C_{I1}^2+C_{I2}\right\} \times \frac{n_I^2}{N_I^2}, \label{app:prelim:sec1:eq3} \\
    \left|Cov(I_{Ii}I_{Ij},I_{Ii}I_{Ij})\right| & \leq \left\{C_{I1}^2+C_{I2}\right\} \times \frac{n_I^2}{N_I^2}, \label{app:prelim:sec1:eq4} \\
    \left|Cov(I_{Ii}I_{Ij},I_{Ii}I_{Ij'})\right| & \leq C'_{I3} \times \frac{n_I^3}{N_I^3} \label{app:prelim:sec1:eq5}, \\
    \left|Cov(I_{Ii}I_{Ij},I_{Ii'}I_{Ij'})\right| & \leq C'_{I4} \times \frac{n_I^3}{N_I^4} \label{app:prelim:sec1:eq6}, \\
    \pi_{Iiji'} & \leq C'_{I5} \frac{n_I^3}{N_I^3}. \label{app:prelim:sec1:eq7}
    \end{align}
\end{lemma}

\subsection{Moment inequalities for the sub-integrals}

\begin{lemma} \label{lem1}
Suppose that assumptions (SS1) and (VA1) hold. Then
    \begin{eqnarray*} \label{app:prelim:eq1}
    \sum_{i=1}^{N_I} \{\tau_{\rho i}\}^4 & \leq & \Phi_1^3 M_1 \times \frac{A^4}{N_I^3}.
    \end{eqnarray*}
\end{lemma}

\noindent By using H\"older's inequality, it also follows from Lemma \ref{lem1} that
    \begin{align}
    \sum_{i=1}^{N_I} \{\tau_{\rho i}\}^2
    & \leq \Phi_1^{\frac{3}{2}} M_1^{\frac{1}{2}} \times \frac{A^2}{N_I}, \label{app:prelim:eq2} \\
    \sum_{i=1}^{N_I} \left|\tau_{\rho i}\right|
    & \leq \Phi_1^{\frac{3}{4}} M_1^{\frac{1}{4}} \times A. \label{app:prelim:eq3}
    \end{align}

\begin{lemma} \label{lem1b}
Suppose that assumptions (FS1), (FS2) and (VA1) hold. Then
    \begin{eqnarray*} \label{app:prelim:eq3b}
     E \left\{\left(\sum_{i \in S_{I}} \frac{\tau_{\rho i}}{\pi_{Ii}} - \tau_{\rho} \right)^4 \right\} & = & O\left(\frac{A^4}{n_I^2} \right).
    \end{eqnarray*}
\end{lemma}

\subsection{Moment inequalities for the second-stage variance}

\begin{lemma} \label{lem2}
Suppose that assumptions (SS1) and (VA1) hold. Then
    \begin{eqnarray*} \label{app:prelim:eq4}
    \sum_{i=1}^{N_I} V_i^2 & \leq & \frac{\Phi_1^3 M_1}{(\lambda_1 c_1)^2} \times \frac{A^4}{N_I^3~n_{II0}^2}.
    \end{eqnarray*}
\end{lemma}

\noindent By using H\"older's inequality, we obtain from (\ref{app:prelim:eq4}) that
    \begin{align} \label{app:prelim:eq5}
    \sum_{i=1}^{N_I} V_i & \leq \frac{\Phi_1^{\frac{3}{2}}~M_1^{\frac{1}{2}}}{\lambda_1 c_1} \times \frac{A^2}{N_I~n_{II0}}.
    \end{align}
If in addition assumption (FS1) holds, then
    \begin{align} \label{app:prelim:eq6}
    \sum_{i=1}^{N_I} \frac{V_i}{\pi_{Ii}} & \leq \frac{\Phi_1^{\frac{3}{2}}~M_1^{\frac{1}{2}}}{\lambda_1 c_1 c_{I1}} \times \frac{A^2}{n_I~n_{II0}}.
    \end{align}

\subsection{Inequalities for higher moments}

\begin{lemma} \label{lem3}
Suppose that assumptions (FS1), (SS1) and (VA1) hold. Then
    \begin{eqnarray}
    \sum_{i=1}^{N_I} V(\Hat{\tau}_{\rho i}^2) & \leq &  \frac{8 M_1}{\lambda_1} \left(\frac{\Phi_1}{c_1} \right)^3 \frac{A^4}{N_I^3~n_{II0}}, \label{app:prelim:eq7} \\
    \sum_{i=1}^{N_I} \frac{V\left(\hat{V}_i\right)}{\pi_{Ii}} & \leq & \frac{16 M_1}{c_{I1}} \left(\frac{\Phi_1}{\lambda_1 c_1} \right)^3 \frac{A^4}{N_I^2 n_I n_{II0}^3}, \label{app:prelim:eq8} \\
     \sum_{i=1}^{N_I} E(\Hat{\tau}_{\rho i}^4) & \leq & 12 \frac{\Phi_1^3 M_1}{\lambda_1^2 c_1^3} \frac{A^4}{N_I^3}. \label{app:prelim:eq8b}
    \end{eqnarray}
\end{lemma}

\begin{lemma} \label{lem4}
Suppose that assumptions (FS1), (FS2), (SS1) and (VA1) hold. Then
    \begin{equation*} \label{app:prelim:eq15}
    E\left(\left[\left\{ \sum_{u_i \in S_{I}} (1-\pi_{Ii}) \frac{\Hat{\tau}_{\rho i}}{\pi_{Ii}} \right\}^2 - \left\{\sum_{i=1}^{N_I} \pi_{Ii}(1-\pi_{Ii}) \frac{\Hat{\tau}_{\rho i}}{\pi_{Ii}} \right\}^2\right]^2 \right) = O \left(\frac{A^4}{n_I^2} \right).
    \end{equation*}
\end{lemma}

\begin{lemma} \label{lem5}
Suppose that assumptions (FS1), (FS2) and (VA1) are respected. Then
    \begin{equation*} \label{app:prelim:eq24}
    E \left\{\left(\Hat{\tau}_{\rho}-{\tau}_{\rho}\right)^4 \right\} = O\left( \frac{A^2}{n_I^2}\right).
    \end{equation*}
\end{lemma}

\begin{lemma} \label{lem5b}
Suppose that assumptions (FS1), (FS2) and (VA1) are respected. Let us denote
    \begin{eqnarray} \label{app:prelim:eq28}
    {\mu}_{\rho} = \frac{{\tau}_{\rho}}{A} & \text{and} & \hat{\mu}_{\rho} = \frac{\hat{\tau}_{\rho}}{\hat{A}},
    \end{eqnarray}
where
    \begin{equation*} \label{app:prelim:eq28b}
    \hat{A} = \sum_{u_i \in S_I} \frac{1}{\pi_{Ii}} \sum_{x \in S_{IIi}} \frac{1}{\pi_{IIi}(x)}.
    \end{equation*}
Then:
    \begin{eqnarray*} \label{app:prelim:eq29}
    E\left\{\left(\hat{\mu}_{\rho}-{\mu}_{\rho} \right)^4 \right\} & = & O\left(\frac{1}{n_I^2} \right).
    \end{eqnarray*}
\end{lemma}

\section{Proof of Proposition \ref{prop1}} \label{app:proof:prop1}

\noindent The variance due to the first-stage is
\begin{align*}
    V_1(\Hat{\tau}_{\rho}) & = \sum_{u_i \in U_I} \sum_{u_j \in U_I} \Delta_{Iij} \frac{\tau_{\rho i}}{\pi_{Ii}} \frac{\tau_{\rho j}}{\pi_{Ij}} \\
    & = \sum_{u_i \in U_I} \frac{1-\pi_{Ii}}{\pi_{Ii}} \{\tau_{\rho i}\}^2+\sum_{u_i \neq u_j \in U_I}  \Delta_{Iij} \frac{\tau_{\rho i}}{\pi_{Ii}} \frac{\tau_{\rho j}}{\pi_{Ij}} \\
    & \leq \frac{1}{\min \pi_{Ii}} \sum_{u_i \in U_I} \{\tau_{\rho i}\}^2 + \frac{\Delta_{I2}}{\{\min \pi_{Ii}\}^2} \left(\sum_{u_i \in U_I} \tau_{\rho i} \right)^2 \quad \text{from line 1 in equation (\ref{sec:assump:eq4})}.
\end{align*}
From equations (\ref{app:prelim:eq2}) and (\ref{app:prelim:eq3}), we have
\begin{align*}
    V_1(\Hat{\tau}_{\rho}) & \leq (\Phi_1)^{\frac{3}{2}} M_1^{\frac{1}{2}}  \times \left[\frac{1}{\min \pi_{Ii}} \times \frac{A^2}{N_I} + \frac{\Delta_{I2}}{\{\min \pi_{Ii}\}^2} \times A^2 \right] \nonumber \\
    & \leq (\Phi_1)^{\frac{3}{2}} M_1^{\frac{1}{2}}  \times \left[\frac{N_I}{c_{I1}~n_I} \times \frac{A^2}{N_I} + C_{I2} \frac{n_I}{N_I^2} \times \frac{N_I^2}{c_{I1}^2~n_I^2} \times A^2 \right] \quad \text{from assumptions (FS1) and (FS2)} \nonumber \\
    & \leq (\Phi_1)^{\frac{3}{2}} M_1^{\frac{1}{2}}  \times \left[\frac{1}{c_{I1}} + \frac{C_{I2}}{c_{I1}^2} \right] \times \frac{A^2}{n_I},
\end{align*}
which leads to (\ref{sec:consist:eq1}). Equation (\ref{sec:consist:eq2}) follows directly from (\ref{app:prelim:eq6}), and equation  (\ref{sec:consist:eq3}) follows directly from (\ref{app:prelim:eq6}).

\section{Proof of Proposition \ref{prop4b}} \label{app:proof:prop4b}

\noindent We begin with equation (\ref{sec:consist:eq10}). From assumption (PI1), we have
    \begin{equation*} \label{proof:prop4b:eq1}
    \hat{\theta}-\theta = g(\hat{\bm{\tau}}_{\bm{\rho}})-g({\bm{\tau}}_{\bm{\rho}}) = g(\hat{\bm{\mu}}_{\bm{\rho}\pi})-g({\bm{\mu}}_{\bm{\rho}}),
    \end{equation*}
where ${\bm{\mu}}_{\bm{\rho}}$ and $\hat{\bm{\mu}}_{\bm{\rho}\pi}$ are defined in equation (\ref{app:prelim:eq28}).
Using assumption (PI2) and the mean-value theorem, there exists some vector $\tilde{\mathbf{a}} \in \mathcal{C}$ whose components lie between those of $\hat{\bm{\mu}}_{\bm{\rho}}$ and ${\bm{\mu}}_{\bm{\rho}}$, and such that
    \begin{eqnarray} \label{proof:prop4b:eq2}
    \hat{\theta}-\theta & = & \left\{\mathbf{g}'(\tilde{\mathbf{a}})\right\}^{\top} \{\hat{\bm{\mu}}_{\bm{\rho}}-{\bm{\mu}}_{\bm{\rho}}\}.
    \end{eqnarray}
From assumption (PI1), the differential $\mathbf{g}'$ is homogeneous of degree $-1$, hence
    \begin{equation} \label{proof:prop4b:eq3}
    \{\mathbf{g}'(\bm{\tau}_{\bm{\rho}})\}^{\top} \{\hat{\bm{\tau}}_{\bm{\rho}}-\bm{\tau}_{\bm{\rho}}\} = \{\mathbf{g}'(\bm{\mu}_{\bm{\rho}})\}^{\top} \{\hat{\bm{\mu}}_{\bm{\rho}}-\bm{\mu}_{\bm{\rho}}\}.
    \end{equation}
From equations (\ref{proof:prop4b:eq2}) and (\ref{proof:prop4b:eq3}), we can write
    \begin{eqnarray*} \label{proof:prop4b:eq4}
    (\hat{\theta}-\theta)-\{\mathbf{g}'(\bm{\tau}_{\bm{\rho}})\}^{\top} \{\hat{\bm{\tau}}_{\bm{\rho}}-\bm{\tau}_{\bm{\rho}}\} & = & \underbrace{\left\{\mathbf{g}'(\tilde{\mathbf{a}})-\mathbf{g}'(\bm{\mu}_{\bm{\rho}})\right\}^{\top} \left\{\hat{\bm{\mu}}_{\bm{\rho}}-\bm{\mu}_{\bm{\rho}}\right\}}_{T_7}.
    \end{eqnarray*}
From the local Lipschitz assumption, we have
    \begin{eqnarray*}
    |T_7| & \leq & \|\mathbf{g}'(\tilde{\mathbf{a}})-\mathbf{g}'(\bm{\mu}_{\bm{\rho}})\| \times \|\hat{\bm{\mu}}_{\bm{\rho}}-\bm{\mu}_{\bm{\rho}}\| \leq K \|\hat{\bm{\mu}}_{\bm{\rho}}-\bm{\mu}_{\bm{\rho}}\|^2,
    \end{eqnarray*}
which gives $E\{|T_7|^2\} \leq K^2 E\left\{\left\|\hat{\bm{\mu}}_{\bm{\rho}}-\bm{\mu}_{\bm{\rho}} \right\|^4 \right\}$. By applying Lemma \ref{lem5} componentwise, we obtain (\ref{sec:consist:eq10}). To prove equation (\ref{sec:consist:eq12}), we first write
    \begin{equation*} \label{proof:prop4b:eq5}
    \hat{V}_{YG}(\hat{\theta})-V_p(\hat{\tau}_{l}) = \underbrace{\hat{V}_{YG}(\hat{\tau}_{l})-V_p(\hat{\tau}_{l})}_{T_8}+\underbrace{\hat{V}_{YG}(\hat{\theta})-\hat{V}_{YG}(\hat{\tau}_{l})}_{T_9}.   \end{equation*}
Since $\mathbf{g}'(\cdot)$ is homogeneous of degree $-1$, we can write
$l(x) = A^{-1} l_0(x)$ with $l_0(x)=\{\mathbf{g}'(\bm{\mu}_{\bm{\rho}})\}^{\top} \bm{\rho}(x)$. From assumption (PI2), there exists some constant $K_c$ such that
    \begin{equation} \label{proof:prop4b:eq6}
    \|\bm{\rho}(x)\| \leq K_c \quad \text{for any } x \in \mcF.
    \end{equation}
Also, the density $l_0(x)$ verifies the first moment condition in assumption (VA1). \\

\noindent We study the two terms in (\ref{proof:prop4b:eq6}) separately. From equation (\ref{sec:consist:eq9}), we directly obtain $E\{|T_8|^2\}=O(n_I^{-3})$ and therefore, $E\{|T_8|\}=O(n^{-3/2})$. Also, we can rewrite $T_9=T_{91}+T_{92}$ with
    \begin{eqnarray} \label{proof:prop4b:eq7}
    T_{91} & = & (\hat{A})^{-2} \times \{\mathbf{g}'(\hat{\bm{\mu}}_{\bm{\rho}})-\mathbf{g}'(\bm{\mu}_{\bm{\rho}})\}^{\top} \times (\mathbf{B}_1+\mathbf{B}_2) \times \{\mathbf{g}'(\hat{\bm{\mu}}_{\bm{\rho}})+\mathbf{g}'(\bm{\mu}_{\bm{\rho}})\}, \nonumber \\
    T_{92} & = & \left\{(\hat{A})^{-2}-A^{-2}\right\} \mathbf{g}'(\bm{\mu}_{\bm{\rho}})^{\top} \times (\mathbf{B}_1+\mathbf{B}_2) \times \mathbf{g}'(\bm{\mu}_{\bm{\rho}}),
    \end{eqnarray}
and where
    \begin{eqnarray*} \label{proof:prop4b:eq8}
    \mathbf{B}_1 & = & -\frac{1}{2} \sum_{u_i \neq u_j \in S_I} \frac{\Delta_{Iij}}{\pi_{Iij}} \left\{\frac{\hat{\bm{\tau}}_{\bm{\rho} i}}{\pi_{Ii}} - \frac{\hat{\bm{\tau}}_{\bm{\rho} j}}{\pi_{Ij}} \right\} \left\{\frac{\hat{\bm{\tau}}_{\bm{\rho} i}}{\pi_{Ii}} - \frac{\hat{\bm{\tau}}_{\bm{\rho} j}}{\pi_{Ij}} \right\}^{\top}, \\
    \mathbf{B}_2 & = & \frac{1}{2} \sum_{u_i \in S_I} \frac{1}{(n_{IIi}-1)\pi_{Ii}} \sum_{x \neq x' \in S_{IIi}} \left\{\frac{\bm{\rho}(x)}{\pi_{IIi}(x)} - \frac{\bm{\rho}(x')}{\pi_{IIi}(x')} \right\} \left\{\frac{\bm{\rho}(x)}{\pi_{IIi}(x)} - \frac{\bm{\rho}(x')}{\pi_{IIi}(x')} \right\}^{\top}. \nonumber
    \end{eqnarray*}
Using assumptions (FS1)-(FS3), we have
    \begin{eqnarray} \label{proof:prop4b:eq9}
    \|\mathbf{B}_1\| & \leq & \frac{C_{I2}}{2c_{I2}(c_{I1})^2} \frac{N_I^2}{n_I^3} \sum_{u_i \neq u_j \in S_I} \left(\left\|\hat{\bm{\tau}}_{\bm{\rho} i}\right\|+\left\|\hat{\bm{\tau}}_{\bm{\rho} j}\right\|\right)^2 \nonumber \\
    & \leq & \frac{2C_{I2}}{c_{I2}(c_{I1})^2} \frac{N_I^2}{n_I^2} \sum_{u_i \in S_I} \left\|\hat{\bm{\tau}}_{\bm{\rho} i}\right\|^2.
    \end{eqnarray}
From equation (\ref{proof:prop4b:eq6}) and using assumption (SS1), we have for any $u_i \in U_I$
$\left\|\hat{\bm{\tau}}_{\bm{\rho} i}\right\| \leq \frac{\Phi_1 K_c A}{c_1 N_I}$,
which together with (\ref{proof:prop4b:eq9}) leads to
    \begin{eqnarray} \label{proof:prop4b:eq11}
    \|\mathbf{B}_1\| & \leq & \frac{2C_{I2}(\Phi_1 K_c)^2}{c_{I2}(c_{I1})^2 (c_1)^2} \frac{A^2}{n_I}.
    \end{eqnarray}
Using assumptions (FS1) and (SS1), we have
    \begin{eqnarray} \label{proof:prop4b:eq12}
    \|\mathbf{B}_2\| & \leq & \frac{\Phi_1^2 A^2}{2 c_{I1}(c_{1})^2 n_{IIi}^2 (n_{IIi}-1) N_I} \sum_{x \neq x' \in S_{IIi}} \left(\left\|\bm{\rho}(x)\right\|+\left\|\bm{\rho}(x')\right\|\right)^2 \nonumber \\
    & \leq & \frac{2 \Phi_1^2 A^2}{c_{I1}(c_{1})^2 n_{IIi}^2 N_I} \sum_{x \in S_{IIi}} \left\|\bm{\rho}(x)\right\|^2 \nonumber \\
    & \leq & \frac{2 \Phi_1^2 (K_c)^2 A^2}{c_{I1}(c_{1})^2 \lambda_1} \times \frac{A^2}{n_{II0} N_I} \quad \text{from (\ref{proof:prop4b:eq6}) and (SS1)}.
    \end{eqnarray}
From equations (\ref{proof:prop4b:eq11}) and (\ref{proof:prop4b:eq12}), there exists some constant $C'$ such that
    \begin{equation} \label{proof:prop4b:eq13}
    \|\mathbf{B}_1+\mathbf{B}_2\| \leq C' \frac{A^2}{n_I}.  \end{equation}

\noindent We first consider the term $T_{91}$ in (\ref{proof:prop4b:eq7}). Since $\mathbf{g}'(\cdot)$ is bounded on $\mathcal{C}$, there exists some constant $C"$ such that
    \begin{eqnarray} \label{proof:prop4b:eq14}
    |T_{91}| & \leq & C" \times (\hat{A})^{-2} \times \|\mathbf{g}'(\hat{\bm{\mu}}_{\bm{\rho}})-\mathbf{g}'(\bm{\mu}_{\bm{\rho}})\| \times \|\mathbf{B}_1+\mathbf{B}_2\| \nonumber \\
    & \leq & C" \times K \times (\hat{A})^{-2} \times \|\hat{\bm{\mu}}_{\bm{\rho}}-\bm{\mu}_{\bm{\rho}}\| \times \|\mathbf{B}_1+\mathbf{B}_2\| \nonumber \\
    & \leq & \frac{C" \times K \times (C_{I1} C_1)^2 \times C'}{(\phi_1)^2 n_I} \times \|\hat{\bm{\mu}}_{\bm{\rho}}-\bm{\mu}_{\bm{\rho}}\|,
    \end{eqnarray}
where the last line in (\ref{proof:prop4b:eq14}) follows from equation (\ref{proof:prop4b:eq13}) and the fact that from assumptions (FS1) and (SS1), we have $\hat{A} \geq \frac{\phi_1 A}{C_{I1} C_1}$. We obtain
    \begin{eqnarray*}
    E\{|T_{91}|\} & \leq & \frac{C" \times K \times (C_{I1} C_1)^2 \times C'}{(\Phi_1)^2 n_I} \times \sqrt{E\{\|\hat{\bm{\mu}}_{\bm{\rho}}-\bm{\mu}_{\bm{\rho}}\|^2\}}. \end{eqnarray*}
Since $E\{\|\hat{\bm{\mu}}_{\bm{\rho}}-\bm{\mu}_{\bm{\rho}}\|^2\}=O(n_I^{-1})$ from Lemma \ref{lem5}, we obtain $E\{|T_{91}|\}=O(n_I^{-3/2})$. Also
    \begin{eqnarray*} \label{proof:prop4b:eq15}
    |T_{92}| & \leq & \frac{\hat{A}+A}{(\hat{A})^2(A)^2} \times |\hat{A}-A| \times \left\|\mathbf{g}'(\bm{\mu}_{\bm{\rho}}) \right\|^2 \times \left\| \mathbf{B}_1+\mathbf{B}_2 \right\| \nonumber \\
    & \leq & \frac{(\Phi_1+c_{I1} c_1) (C_{I1} C_1)^2}{\Phi_1^2 (c_{I1} c_1) A^3} \times \left\|\mathbf{g}'(\bm{\mu}_{\bm{\rho}}) \right\|^2 \times \left\| \mathbf{B}_1+\mathbf{B}_2 \right\| \times |\hat{A}-A|,
    \end{eqnarray*}
which leads to
    \begin{equation*} \label{proof:prop4b:eq16}
    E\left\{|T_{92}|\right\} \leq \frac{(\Phi_1+c_{I1} c_1) (C_{I1} C_1)^2}{\Phi_1^2 (c_{I1} c_1) A} \times \frac{(C")^2 C'}{n_I} \times \sqrt{E\left\{(\hat{A}-A)^2\right\}}.
    \end{equation*}
From equation (\ref{sec:consist:eq4}) applied with $\rho(x)=1$, we have $E\left\{(\hat{A}-A)^2\right\}=O(A^2 n_I^{-1})$. This leads to $E\{|T_{92}|\}=O(n_I^{-3/2})$, which completes the proof for equation (\ref{sec:consist:eq12}). The proof for equation (\ref{sec:consist:eq11}) is similar.

\section{Proof of Proposition \ref{prop5}} \label{app:proof:prop5}

\noindent We first prove equation (\ref{sec:large:entropy:eq3}). We simplify the notation as $\hat{V}_{HAJ,A} \equiv \hat{V}_{HA}$, $V_1(\Hat{\tau}_{\rho}) \equiv V_1$ and $V_2(\Hat{\tau}_{\rho}) \equiv V_2$. We also use the notation
    \begin{align*}
      \tilde{V}_{HA} & = \sum_{i=1}^{N_I} \pi_{Ii}(1-\pi_{Ii}) \left(\frac{\Hat{\tau}_{\rho i}}{\pi_{Ii}} - \hat{R} \right)^2 \quad \text{with} \quad \hat{R} = d_I^{-1} \sum_{i=1}^{N_I} \pi_{Ii}(1-\pi_{Ii}) \frac{\Hat{\tau}_{\rho i}}{\pi_{Ii}}, \\
      \tilde{V}_{HT,A} & = \sum_{i,j=1}^{N_I} \frac{\Delta_{Iij}}{\pi_{Ii}\pi_{Ij}} \Hat{\tau}_{\rho i} \Hat{\tau}_{\rho j} = E \left[ \left. \hat{V}_{HT,A}(\Hat{\tau}_{\rho}) \right| S_{II}\right],
    \end{align*}
and $S_{II} = \cup_{i=1}^{N_I} S_i$ for the union of the second-stage samples. We write
    \begin{equation} \label{app:proof:prop5:eq3}
    \hat{V}_{HAJ,A}-V_1-V_2 = T_{A1}-T_{A2}-T_{A3}+T_{A4}+T_{A5},
    \end{equation}
with
    \begin{align*}
      T_{A1} & = \sum_{i=1}^{N_I} (I_{Ii}-\pi_{Ii})(1-\pi_{Ii}) \left(\frac{\Hat{\tau}_{\rho i}}{\pi_{Ii}} \right)^2 , \\
      T_{A2} & = (\hat{d}_{rI}^{-1}-d_{I}^{-1}) \left\{ \sum_{i=1}^{N_I} \pi_{Ii}(1-\pi_{Ii}) \frac{\Hat{\tau}_{\rho i}}{\pi_{Ii}} \right\}, \\
      T_{A3} & = d_{I}^{-1} \left[\left\{ \sum_{i \in S_{rI}} (1-\pi_{Ii}) \frac{\Hat{\tau}_{\rho i}}{\pi_{Ii}} \right\}^2 - \left\{\sum_{i=1}^{N_I} \pi_{Ii}(1-\pi_{Ii}) \frac{\Hat{\tau}_{\rho i}}{\pi_{Ii}} \right\}^2\right],
    \end{align*}
and with $T_{A4} =\tilde{V}_{HT,A} -V_1-V_2$ and $T_{A5} =\tilde{V}_{HA}-\tilde{V}_{HT,A}$. For the first term in (\ref{app:proof:prop5:eq3}), we have
    \begin{eqnarray*}
      E(T_{A1}^2) & \leq & E\left[ \sum_{i=1}^{N_I} (I_{Ii}-\pi_{Ii})(1-\pi_{Ii}) \left(\frac{\Hat{\tau}_{\rho i}}{\pi_{Ii}} \right)^2\right]^2 \nonumber \\
                  & = & EV\left[ \left.\sum_{i=1}^{N_I} (I_{Ii}-\pi_{Ii})(1-\pi_{Ii}) \left(\frac{\Hat{\tau}_{\rho i}}{\pi_{Ii}} \right)^2 \right| S_{II}\right] \nonumber \\
                  & \leq & E\left[ \sum_{i=1}^{N_I} \frac{(1-\pi_{Ii})^3}{\pi_{Ii}^3} (\Hat{\tau}_{\rho i})^4 \right] \leq \sum_{i=1}^{N_I} \frac{E(\Hat{\tau}_{\rho i})^4}{\pi_{Ii}^3}, \label{app:proof:prop5:eq4}
    \end{eqnarray*}
and from assumption (FS1) and equation (\ref{app:prelim:eq8b}) in Lemma \ref{lem3}, we obtain
    \begin{eqnarray} \label{app:proof:prop5:eq5}
      E(T_{A1}^2) & = & O\left(\frac{A^4}{n_I^3}\right).
    \end{eqnarray}

\noindent For the second term in (\ref{app:proof:prop5:eq3}), we obtain from Lemma \ref{lem0}
    \begin{eqnarray*}
      E(T_{A2}^2) & \leq & \frac{4}{(c_{I1})^4 (1-f_{I0})^4 n_I^4} E(\hat{d}_{rI}-d_I)^2 \times E \left\{ \sum_{i=1}^{N_I} (1-\pi_{Ii}) \Hat{\tau}_{\rho i} \right\}^4 \nonumber \\
      & \leq & \frac{4}{(c_{I1})^4 (1-f_{I0})^4 n_I^3} \times N_I^3 \sum_{i=1}^{N_I} E(\Hat{\tau}_{\rho i})^4, \label{app:proof:prop5:eq6}
    \end{eqnarray*}
and from equation (\ref{app:prelim:eq8b}) in Lemma \ref{lem3}, we obtain
    \begin{eqnarray} \label{app:proof:prop5:eq7}
      E(T_{A2}^2) & = & O\left(\frac{A^4}{n_I^3}\right).
    \end{eqnarray}

\noindent \noindent For the third term in (\ref{app:proof:prop5:eq3}), we obtain from Lemmas \ref{lem0} and \ref{lem3} that
    \begin{eqnarray} \label{app:proof:prop5:eq8}
      E(T_{A3}^2) & = & O\left(\frac{A^4}{n_I^4}\right).
    \end{eqnarray}

\noindent \noindent For the fourth term in (\ref{app:proof:prop5:eq3}), we have
    \begin{equation*} \label{app:proof:prop5:eq9}
      E(T_{A4}^2) = V\left\{\tilde{V}_{HT,A}\right\} \leq 2 V\left\{\sum_{i=1}^{N_I} \frac{1-\pi_{Ii}}{\pi_{Ii}} (\Hat{\tau}_{\rho i})^2 \right\} + 2 V\left\{\sum_{i \neq j =1}^{N_I} \frac{\Delta_{Iij}}{\pi_{Ii}\pi_{Ij}} \Hat{\tau}_{\rho i} \Hat{\tau}_{\rho j} \right\}.
    \end{equation*}
We have
    \begin{eqnarray}  \label{app:proof:prop5:eq10}
      V\left\{\sum_{i=1}^{N_I} \frac{1-\pi_{Ii}}{\pi_{Ii}} (\Hat{\tau}_{\rho i})^2 \right\} & \leq & \sum_{i=1}^{N_I} \frac{V\{(\Hat{\tau}_{\rho i})^2\}}{\pi_{Ii}^2} =O\left(\frac{A^4}{n_I^3}\right),
    \end{eqnarray}
where the order of magnitude in (\ref{app:proof:prop5:eq10}) follows from assumption (FS1) and from equation (\ref{app:prelim:eq7}) in Lemma \ref{lem3}. We also have
    \begin{align} \label{app:proof:prop5:eq11}
      V\left\{\sum_{i \neq j =1}^{N_I} \frac{\Delta_{Iij}}{\pi_{Ii}\pi_{Ij}} \Hat{\tau}_{\rho i} \Hat{\tau}_{\rho j} \right\} & = 2 \sum_{i \neq j \neq j'=1}^{N_I} \frac{\Delta_{Iij}}{\pi_{Ii}\pi_{Ij}} \frac{\Delta_{Iij'}}{\pi_{Ii}\pi_{Ij'}} Cov\left\{\Hat{\tau}_{\rho i}\Hat{\tau}_{\rho j},\Hat{\tau}_{\rho i}\Hat{\tau}_{\rho j'}\right\} + \sum_{i \neq j =1}^{N_I} \left(\frac{\Delta_{Iij}}{\pi_{Ii}\pi_{Ij}}\right)^2 V(\Hat{\tau}_{\rho i}\Hat{\tau}_{\rho j}) \nonumber \\
      & = 2 \sum_{i \neq j \neq j'=1}^{N_I} \frac{\Delta_{Iij}}{\pi_{Ii}\pi_{Ij}} \frac{\Delta_{Iij'}}{\pi_{Ii}\pi_{Ij'}} \left\{{\tau}_{\rho j} {\tau}_{\rho j'} V_i\right\} + \sum_{i \neq j =1}^{N_I} \left(\frac{\Delta_{Iij}}{\pi_{Ii}\pi_{Ij}}\right)^2 \left\{{\tau}_{\rho j}^2 V_i + {\tau}_{\rho i}^2 V_j + V_i V_j \right\} \nonumber \\
      & \leq \frac{C_{I2}^2}{c_{I1}^4} \times \frac{1}{n_I^2} \times \left[ 2 \left\{\sum_{i=1}^{N_I} {\tau}_{\rho i} \right\}^2 \sum_{j=1}^{N_I} V_j + 2 \sum_{i=1}^{N_I} {\tau}_{\rho i}^2 \sum_{j=1}^{N_I} V_j + \left\{\sum_{j=1}^{N_I} V_j\right\}^2 \right] \nonumber \\
      & = O\left(\frac{A^4}{n_I^3}\right),
    \end{align}
where the order of magnitude in (\ref{app:proof:prop5:eq11}) follows from equations (\ref{app:prelim:eq2}), (\ref{app:prelim:eq3}) and (\ref{app:prelim:eq5}). This leads to
    \begin{eqnarray} \label{app:proof:prop5:eq12}
      E(T_{A4}^2) = V\left\{\tilde{V}_{HT,A}\right\} =O\left(\frac{A^4}{n_I^3}\right).
    \end{eqnarray}

\noindent For the fifth term in (\ref{app:proof:prop5:eq3}), we have from \citet[][Theorem 5.2]{Haj64} $\tilde{V}_{HA} = \tilde{V}_{HT,A}\{1+o(1)\}$, which leads to $E(T_{A5}^2) \leq  E\left\{(\tilde{V}_{HT,A})^2\right\} \times o(1)$. We have
    \begin{eqnarray} \label{app:proof:prop5:eq15}
      E\left\{(\tilde{V}_{HT,A})^2\right\} & = & V\left\{\tilde{V}_{HT,A}\right\} + \left[V_1+V_2\right]^2 = O\left(\frac{A^4}{n_I^2} \right),
    \end{eqnarray}
where the order of magnitude in (\ref{app:proof:prop5:eq15}) follows from equation (\ref{app:proof:prop5:eq12}) and Proposition \ref{prop1}. This leads to $E(T_{A5}^2) = o\left(\frac{A^4}{n_I^2} \right)$. This completes the proof.

\noindent Now, we consider equation (\ref{sec:large:entropy:eq4}). From Section 3.2 in the supplementary material of \citet{chauvet2020inference}, it is sufficient to prove that
    \begin{eqnarray}
    \frac{\sum_{i=1}^{N_I} \frac{E(\Hat{\tau}_{\rho i}-{\tau}_{\rho i})^4}{\pi_{Ii}^3}}{\left\{V(\Hat{\tau}_{\rho})\right\}^2} & = & o(1), \label{app:proof:prop5:eq17} \\
    \frac{\sum_{i=1}^{N_I} \frac{({\tau}_{\rho i}-R \pi_{Ii})^4}{\pi_{Ii}^3}}{\left\{V(\Hat{\tau}_{\rho})\right\}^2} & = & o(1), \label{app:proof:prop5:eq18}
    \end{eqnarray}
where $R=d_I^{-1} \sum_{i=1}^{N_I} \pi_{Ii}(1-\pi_{Ii}) \frac{\tau_{\rho i}}{\pi_{Ii}}$. From assumption (FS1) and equation (\ref{app:prelim:eq8b}) in Lemma \ref{lem3}, we have
    \begin{eqnarray} \label{app:proof:prop5:eq19}
    \sum_{i=1}^{N_I} \frac{E(\Hat{\tau}_{\rho i}-{\tau}_{\rho i})^4}{\pi_{Ii}^3} & = & O\left(\frac{A^4}{n_I^3}\right),
    \end{eqnarray}
and equation (\ref{app:proof:prop5:eq17}) follows from (\ref{app:proof:prop5:eq19}) and assumption (VA2). From assumption (FS1) and from Lemma \ref{lem1}, we have
    \begin{eqnarray} \label{app:proof:prop5:eq20}
    \sum_{i=1}^{N_I} \frac{({\tau}_{\rho i}-R \pi_{Ii})^4}{\pi_{Ii}^3} & = & O\left(\frac{A^4}{n_I^3}\right),
    \end{eqnarray}
and equation (\ref{app:proof:prop5:eq18}) follows from (\ref{app:proof:prop5:eq20}) and assumption (VA2). Finally, we consider (\ref{sec:large:entropy:eq5}). It follows directly from (\ref{sec:consist:eq6}) and (\ref{sec:large:entropy:eq3}) and assumption (VA2).

\section{Proof of Proposition \ref{prop6}} \label{app:proof:prop6}

\noindent The proof is based in the coupling procedure given in \citet[][Algorithm 1]{chauvet2020inference}. Using the same steps as in Section 4 of the supplementary material of \citet[][Algorithm 1]{chauvet2020inference}, it is sufficient to prove that
    \begin{align}
    E\left\{\left(\sum_{i \in S_{rI}} \frac{\tau_{\rho i}}{\pi_{Ii}} - \tau_{\rho} \right)^4 \right\} & = O\left(\frac{A^4}{n_I^2} \right), \label{proof:prop6:eq1} \\
    E \left\{\left(\sum_{i \in S_{rI}} \frac{V_i}{\pi_{Ii}^2}\right)^2\right\} & = O\left(\frac{A^4}{n_I^2} \right). \label{proof:prop6:eq2}
    \end{align}
Equation (\ref{proof:prop6:eq1}) is established in Lemma \ref{lem1b}. For equation (\ref{proof:prop6:eq2}), we have by the Cauchy-Schwartz inequality
    \begin{equation} \label{proof:prop6:eq3}
    E \left\{\left(\sum_{i \in S_{rI}} \frac{V_i}{\pi_{Ii}^2}\right)^2  = O\left(\frac{A^4}{n_I^2} \right) \right\} \leq n_I \sum_{i=1}^{N_I} \frac{V_i}{\pi_{Ii}^3} = O\left(\frac{A^4}{n_I^2 n_{II0}} \right),
    \end{equation}
where the order of magnitude in (\ref{proof:prop6:eq3}) follows from assumption (FS2) and Lemma \ref{lem2}.

\section{Proof of Proposition \ref{prop7}} \label{app:proof:prop7}

\noindent Write $\Delta=\bm{g}(\hat{\bm{\tau}}_{\bm{\rho}})-\bm{g}(\bm{\tau}_{\bm{\rho}})-\{\bm{g}'(\bm{\tau}_{\bm{\rho}})\}^{\top} \{\hat{\bm{\tau}}_{\bm{\rho}}-\bm{\tau}_{\bm{\rho}}\}$. We have the identity
    \begin{equation} \label{app:proof:prop7:eq1}
    \frac{\hat{\theta}-\theta}{\sqrt{V(\hat{\theta})}} = \frac{\hat{\tau}_l - \tau_l}{\sqrt{V(\hat{\tau}_l)}} \times \sqrt{\frac{V(\hat{\tau}_l)}{V(\hat{\theta})}} + \frac{\Delta}{\sqrt{V(\hat{\theta})}}.
    \end{equation}
From equations (\ref{sec:consist:eq10}) and (\ref{prop7:eq1}), the second term in the right-hand side of (\ref{app:proof:prop7:eq1}) is $o_p(1)$ and $V(\hat{\tau}_l)^{-1}V(\hat{\theta}) \rightarrow_{Pr} 1$. From equation (\ref{sec:large:entropy:eq7}) and Slutsky's theorem, we obtain (\ref{prop7:eq2}). To prove equation (\ref{prop7:eq3}), we obtain from (\ref{sec:large:entropy:eq8}) that
    \begin{equation*} 
      E\left[ A^{-2} n_I \left| \hat{V}_{HAJ}(\hat{\tau}_l)-V(\Hat{\theta})\right| \right] = o(1).
    \end{equation*}
It is therefore sufficient to prove that $E\left[ A^{-2} n_I \left| \hat{V}_{HAJ}(\hat{\theta})-\hat{V}_{HAJ}(\hat{\tau}_l)\right| \right] = o(1)$. The proof is similar to that of equation (\ref{sec:consist:eq12}).

\bibliographystyle{apalike}

\clearpage
\section*{Author information}

\noindent
Guillaume Chauvet, corresponding author\\
Univ Rennes, Ensai, CNRS, CREST – UMR 9194, F-35000 Rennes, France\\
\href{mailto:guillaume.chauvet@ensai.fr}{guillaume.chauvet@ensai.fr}\\[0.8em]

\noindent
Olivier Bouriaud\\
\c{S}tefan cel Mare University of Suceava, Romania\\
Universit{\'e} de Lorraine, G{\'e}odata Paris, IGN, LIF, F-54000, Nancy, France\\
Universit{\'e} Gustave Eiffel, G{\'e}odata Paris, IGN, Laboratoire d'Inventaire Forestier, F-54000, Nancy, France \\
\href{mailto:obouriaud.lif@gmail.com}{obouriaud.lif@gmail.com}\\[0.8em]

\noindent
Trinh H.K. Duong\\
Universit{\'e} de Lorraine, G{\'e}odata Paris, IGN, LIF, F-54000, Nancy, France\\
Universit{\'e} Gustave Eiffel, G{\'e}odata Paris, IGN, Laboratoire d'Inventaire Forestier, F-54000, Nancy, France \\
\href{mailto:kimtrinh0311@gmail.com}{kimtrinh0311@gmail.com}

\end{document}